\documentclass[10pt]{article}
\newcommand{\heuteIst}{November 11, 2002 }
\usepackage{amsmath,amsthm,amsfonts,latexsym,amscd,amssymb}
\usepackage[dvips]{hyperref}

\scrollmode

\swapnumbers
\theoremstyle{plain}
\newtheorem{theorem}{Theorem}[section]
\newtheorem{lemma}[theorem]{Lemma}
\newtheorem{corollary}[theorem]{Corollary}
\newtheorem{proposition}[theorem]{Proposition}

\newtheorem{definition}[theorem]{Definition}
\theoremstyle{definition}
\newtheorem{remark}[theorem]{Remark}
\newtheorem{question}[theorem]{Question}

\newcommand{\reals}{\mathbb{R}}

\newcommand{\complexs}{\mathbb{C}}

\newcommand{\naturals}{\mathbb{N}}
\newcommand{\Z}{\mathbb{Z}}
\newcommand{\integers}{\mathbb{Z}}
\newcommand{\rationals}{\mathbb{Q}}
\newcommand{\Q}{\mathbb{Q}}

\DeclareMathOperator{\id}{id}
\newcommand{\boundary}[1]{\partial#1}

\newcommand{\abs}[1]{\left\lvert#1\right\rvert} 
\newcommand{\norm}[1]{\lVert#1\rVert}

\newcommand{\tensor}{\otimes}
\newcommand{\into}{\hookrightarrow}

\newcommand{\iso}{\cong}
\newcommand{\disjointunion}{\amalg}

\newcommand{\NeumannN}{\mathcal{N}}

\newcommand{\innerprod}[1]{\langle #1 \rangle}
\DeclareMathOperator{\supp}{supp}   
\DeclareMathOperator{\im}{im}      

\DeclareMathOperator{\Hom}{Hom}    

\DeclareMathOperator{\tr}{tr}

\DeclareMathOperator{\coker}{coker}

\DeclareMathOperator{\ind}{ind}
\DeclareMathOperator{\sign}{sign}

\DeclareMathOperator{\smax}{max}  
\DeclareMathOperator{\san}{an}  
\DeclareMathOperator{\schain}{chain}  
\DeclareMathOperator{\sforms}{forms}  

\newcommand{\forget}[1]{}

{\catcode`@=11\global\let\c@equation=\c@theorem}

\allowdisplaybreaks[2]

\newcommand{\squarematrix}[4]                
{                                            
\begin{pmatrix} #1 & #2 \\ #3 &
#4
\end{pmatrix}
}

\begin{document}

\date{Last edited \heuteIst or later --- last compiled: \today}

\title{Various $L^2$-signatures and a topological $L^2$-signature theorem}
\author{ Wolfgang
L{\"u}ck\thanks{\noindent email:
lueck@math.uni-muenster.de\protect\\
www: ~http://www.math.uni-muenster.de/u/lueck/org/staff/lueck/\protect}\\
Fachbereich Mathematik\\ Universit{\"a}t M{\"u}nster\\
Einsteinstr.~62\\         48149 M{\"u}nster \and
Thomas Schick\thanks{\noindent e-mail:
schick@uni-math.gwdg.de\protect\\
www:~http://www.uni-math.gwdg.de/schick\protect\\
Research partially carried out during a stay at Penn State university
funded by the DAAD}\\ Fakult{\"a}t f{\"u}r Mathematik\\ Universit{\"a}t
G{\"o}ttingen\\ Bunsenstrasse 3\\
37073 G{\"o}ttingen\\ Germany}
\maketitle

\typeout{--------------------  Abstract  --------------------}

\begin{abstract}
  For a normal covering over a closed oriented topological manifold 
  we give a proof of the
  $L^2$-signature theorem with twisted coefficients,
  using Lipschitz structures and the Lipschitz signature operator
  introduced by Teleman. We also prove that the $L$-theory isomorphism
  conjecture as well as the $C^*_{\smax}$-version of the
  Baum-Connes conjecture imply the $L^2$-signature theorem for 
  a normal covering over a  Poincar{\'e} space, provided that the 
group of deck transformations is torsion-free.

  We discuss the various possible
  definitions of $L^2$-signatures (using the signature operator, using 
  the cap product of differential forms, using a cap product in
  cellular $L^2$-cohomology,~\ldots) in this situation, and prove that
  they all coincide.

\end{abstract}

\noindent
Key words: $L^2$-signature, signature, Lipschitz manifolds,
$L^2$-signature theorem.
\\
2000 mathematics subject classification: 57P10, 
57N65, 
58G10. 

\typeout{--------------------   Introduction  --------------------}

\setcounter{section}{-1}
\section{Introduction}
\label{sec:intro}

Atiyah's celebrated $L^2$-index theorem \cite{Atiyah(1976)} implies 
that the index of the signature operator of a closed oriented
smooth manifold $M$ with Riemannian metric coincides with the $L^2$-index of the signature
operator on any normal covering space of $M$. In particular, the
signature and the $L^2$-signature for closed oriented smooth manifolds
coincide. The (various) definitions of $L^2$-signatures are explained
in Section \ref{sec:L2sig}.

The signature is of course also defined for closed oriented topological manifolds
and, as long as there is a Lipschitz structure, there is even a
signature operator whose index is the signature. In the first part of
this paper, we address the question how to generalize the
$L^2$-signature theorem to closed oriented topological manifolds.

Such an $L^2$-signature theorem for closed oriented topological manifolds does not
seem to be in the literature. We give a 
proof along the lines of Atiyah's proof \cite{Atiyah(1976)} of the smooth
$L^2$-index theorem.
  \begin{theorem}
    Let $M$ be a closed connected oriented $4n$-dimensional Lipschitz manifold with normal covering 
    $\overline{M} \to M$. Let $D_V$ be the Lipschitz 
   signature operator twisted with a
    Lipschitz bundle $V$ and $\overline{D_V}$ its lift to $\overline{M}$. Then
    \begin{equation*}
      \ind(D_V) = \ind_{\NeumannN\Gamma}(\overline{D_V}).
    \end{equation*}
  \end{theorem}
An immediate consequence is
(using Sullivan's theorem that a Lipschitz structure exists on every
topological manifold of dimension $\ne 4$)
  \begin{theorem} \label{sign^{(2)}(bar M) = sign(M)}
    Let $M$ be a closed connected oriented topological manifold of dimension $4n$ with normal covering
    $\overline{M} \to M$. Then
    \begin{equation*}
      \sign^{(2)}(\overline{M}) = \sign(M).
    \end{equation*}
  \end{theorem}
Theorem \ref{sign^{(2)}(bar M) = sign(M)}
 also follows from the $L^2$-signature theorem for
closed oriented smooth manifolds and the fact that the forgetful map
$\Omega_*(B\Gamma) \to \Omega^{\text{top}}_* (B\Gamma)$
from the smooth bordism group over
$B\Gamma$ to the topological one is rationally an isomorphism (compare
 Remark \ref{rem:Shmuels_proof} and
 the discussion after \cite[Theorem 1.6]{Teleman(1984)}), as was
 pointed out to us by Shmuel Weinberger.
 Note that the $L^2$-signature theorem implies in particular that the
 signature is multiplicative under finite coverings. This multiplicativity was
 proved for closed oriented topological manifolds in \cite[Theorem
 8]{Schafer(1970)}.

For more general Poincar{\'e} duality spaces, which are not manifolds,  such a multiplicativity
result does not hold \cite[Example 22.28]{Ranicki(1992)}, \cite[Corollary
5.4.1]{Wall(1967)}). It fails also if $M$ is a compact oriented smooth manifold with nonempty boundary
(compare \cite[Proposition 2.12]{Atiyah-Patodi-Singer(1975b)}
together with the Atiyah-Patodi-Singer index
theorem \cite[Theorem 4.14]{Atiyah-Patodi-Singer(1975a)}). 

This implies in particular that the $L^2$-index theorem can not hold
in the stated form in the greatest imaginable
generality. In Section \ref{sec:poincarecomp}, we discuss to which extent the $L^2$-signature
theorem does extend to  Poincar{\'e} spaces $X =
(X,\emptyset)$, and show that it is implied by the $L$-theory
isomorphism conjecture or by the $C_{\smax}^*$-version of
the Baum-Connes conjecture, provided the covering group $\Gamma$ is
torsion-free. More precisely, we prove the following theorem in
\ref{pcl2sig}:
\begin{theorem}
  Let $X$ be a $4n$-dimensional 
  Poincar{\'e} space over $\rationals$ (see Definition \ref{def: Poincare pair}).
  Let $\overline{X}\to X$ be a normal
  covering with \emph{torsion-free} covering group $\Gamma$. Assume
  that the (Baum-Connes)
  index map for the maximal group $C^*$-algebra
  \begin{equation*}
     \ind\colon K_0(B\Gamma) \to K_0(C^*_{\smax}\Gamma)
   \end{equation*}
   or the $L$-theory assembly map
 \begin{equation*}
     A\colon H_{4n}(B\Gamma;{\mathbb L}^{\langle - \infty \rangle}_{\bullet}) \to 
   L_0^{\langle - \infty \rangle} (\integers\Gamma)
   \end{equation*}
   is rationally surjective.
   Then
   \begin{equation*}
     \sign^{(2)}(\overline{X}) = \sign(X).
   \end{equation*}
\end{theorem}

In a companion paper \cite{Lueck-Schick(2001)} we show that, without
any further assumption, multiplicativity of $L^2$-signatures under
coverings holds
``approximately'' in the following sense:
\begin{theorem} \cite[Theorem 0.1]{Lueck-Schick(2001)}
Let $(X,Y)$ be a  $4n$-dimensional Poincar{\'e} pair  over $\rationals$.
Suppose that there is a nested sequence of normal subgroups of finite index
$\Gamma \supseteq \Gamma_1 \supseteq \Gamma_2 \supseteq \ldots$
such that the intersection of the $\Gamma_k$-s is trivial. Let $(X_k,Y_k) \to (X,Y)$ be the 
finite covering of $X$ associated to $\Gamma_K \subseteq \Gamma$.
Then the sequence
$\left(\sign(X_k,Y_k)/[\Gamma : \Gamma_k]\right)_{k \ge 1}$
converges and
  \begin{equation*}
   \lim_{k\to \infty}    \frac{\sign(X_k,Y_k)}{[\Gamma : \Gamma_k]} =
   \sign^{(2)}(\overline{X},\overline{Y}).
  \end{equation*}
\end{theorem}
(In \cite{Lueck-Schick(2001)}, we also prove a similar approximation
result for amenable exhaustions).

In the last part of the present paper, we check that the various versions
of $L^2$-signatures, e.g.~given in terms of intersection pairings, the
index of the signature operator, or as trace of an index element in
the $K$-theory of certain $C^*$-algebras, all coincide whenever the
definitions make sense. In the rest of the paper, and also in
\cite{Lueck-Schick(2001)}, we already freely jump
between the different interpretations.

This comparison is even interesting for smooth manifolds, in particular for
smooth manifolds with boundary. In this case, the
$L^2$-signature is
defined in terms of the intersection pairing on $L^2$-homology. In
Theorem \ref{the: L2indtheorem} we give a proof that this coincides
with the answer predicted by the $L^2$-index theorem \cite[Theorem
1.1]{Ramachandran(1993)}. Note that we deliberately write ``answer
predicted by the $L^2$-index theorem'' and not ``index of the
signature operator'', because before adding a certain
well defined correction term (compare
\cite{Atiyah-Patodi-Singer(1975c)}) one can not expect to obtain the
signature. The paper \cite{Ramachandran(1993)} only deals with the
$L^2$-index
of certain operators. The homological interpretation does not seem to
have been checked in the literature.

\textbf{Organization of the paper:}
In Section \ref{sec:topo} we prove the $L^2$-signature theorem for closed
topological manifolds.

In Section \ref{sec:poincarecomp} we address
the question, for which Poincar{\'e} spaces an $L^2$-signature
theorem holds.

In Section \ref{sec:L2sig}, we compare the different definitions of
$L^2$-signatures, and show that they all coincide.


\section{$L^2$-signature theorem for topological manifolds}\label{sec:topo}

We prove the $L^2$-signature theorem for closed oriented Lipschitz
manifolds. This does prove
the theorem for arbitrary oriented topological manifolds because Sullivan
constructs  in dimensions
$\ge 5$ a (unique) Lipschitz structure on every topological
manifold \cite{Sullivan(1979)}, and taking the product with
$\complexs P^8$ if necessary (which does
change neither the signature nor the $L^2$-signature (compare
Proposition \ref{sigmult})), we may assume
that the dimension of our manifold is sufficiently high.
Note that we need only the existence, but not the uniqueness of this
Lipschitz structure.

Now suppose that $M$ is a closed connected oriented Lipschitz manifold
of dimension $4n$ with a Lipschitz
metric $g$ and with fundamental group $\Gamma$. Let $V$ be a finite
dimensional Lipschitz Hermitian vector bundle over $M$ with a (not necessarily
flat) Lipschitz connection. Teleman \cite{Teleman(1983)}
constructs then a twisted signature operator $D_V$  (whose index is
the topological signature of $M$ if $V$ is a trivial flat line bundle).
For basics about Lipschitz manifolds, Lipschitz bundles and Lipschitz
operators compare  \cite[Section
1--6]{Teleman(1983)}, \cite[Section 2]{Hilsum(1985)}, \cite[Section
1]{Hilsum(1989)}.
 The
Lipschitz structure, the metric, the bundle, and the signature
operator all can be
lifted to $\overline{M}$ and then in particular $\ind_{\NeumannN\Gamma}(\overline{D})=
\sign^{(2)}(\overline{M})$ (if again $V$ is a trivial line bundle). The
task is now to compare $\ind(D_V)$ and
$\ind_{\NeumannN\Gamma}(\overline{D_V})$, which in the smooth case is
done in Atiyah's
$L^2$-index theorem \cite[(1.1)]{Atiyah(1976)}.

The subscript $\NeumannN\Gamma$ refers to the group von Neumann algebra.
Basics about  $\NeumannN\Gamma$, Hilbert $\NeumannN (\Gamma)$-modules, the standard trace 
$\tr_{\NeumannN\Gamma}$ and the von Neumann dimension
$\dim_{\NeumannN\Gamma}$ as used
in this paper can be found e.g.~ in \cite[Section 1 and 2]{Lott-Lueck(1995)},
\cite[Section 1.1]{Lueck(2002)}.

On Lipschitz manifolds, no
pseudo-differential calculus in the usual sense exists.
However, one has the following properties of the
twisted signature operator which are essential for Atiyah's proof:
\begin{theorem}\label{opprop}
  \begin{enumerate}
  \item \label{selfadjoint}
    (The closure of) $D_V$ is an unbounded
    selfadjoint operator. The same is true for $\overline{D_V}$.
  \item \label{unitspeed} $D_V$ and $\overline{D_V}$ have unit propagation
    speed, i.e.
    \begin{equation*}
\supp(e^{itD_V})
    \subseteq \{(x,y);\; (x,y)\in M\times M\text{ and }d(x,y)\le t\}
  \end{equation*}
  and correspondingly for $\overline{D_V}$.
  \item \label{summability} $(i+D)^{-1}$ is $(\dim(M)+1)$-summable.
\end{enumerate}
\end{theorem}
\begin{proof}
These results are established by Hilsum for the untwisted
Lipschitz-signature operator on a complete oriented Lipschitz
manifold. Hilsum uses specific properties of the
untwisted operator so that his proofs can not directly be applied. We
reduce the twisted case to the untwisted case in the following way:
we embed $V$ as Hermitian vector bundle in an $N$-dimensional trivial bundle
$\underline{N}$ with complement $W$. Choose a connection on $W$. We
then have on $\underline{N}$ the trivial connection and the direct sum
of these two connections. Correspondingly, we get two twisted
signature operators $D_N$ (which of course is the $N$-fold direct sum
of the untwisted signature operator) and $D_V\oplus D_W$. A
calculation in local coordinates shows that
\begin{equation*}
  D_V\oplus D_W = D_N + A
\end{equation*}
where $A$ is a bundle homomorphism with bounded measurable and selfadjoint 
coefficients,
therefore is a bounded selfadjoint operator on
$L^2(\Omega^*(M,\underline{N}))$ (for all this compare \cite[Section 6
and 7]{Teleman(1983)}. In fact, in \cite{Teleman(1983)} this is used as the definition of the
twisted signature operator). Lifting this gives the corresponding
splitting
\begin{equation*}
  \overline{D_V}\oplus \overline{D_W} = \overline{D_N} + \overline{A}.
\end{equation*}

By \cite[Corollaire 1.8]{Hilsum(1989)}
the untwisted operators $D_N$
and $\overline{D_N}$ and then also
$D_N + A$ and $\overline{D_N} + \overline{A}$ are selfadjoint, therefore
the same is true for the direct summands $D_V$ and $\overline{D_V}$.

For summability we have to find a relation between $(D_N+i)^{-1}$
(which is $(\dim(M)+1)$-summable by \cite[Proposition 5.6]{Hilsum(1985)}) 
and
$(D_N+A+i)^{-1}$ (these are bounded operators because of
self-adjointness). We compute
\begin{equation*}
    (D_N+A+i)^{-1} =   (D_N+i)^{-1}(1+A(D_N+i)^{-1})^{-1}.
  \end{equation*}
  Since the space of $(\dim(M)+1)$-summable operators is an ideal in the space
  of bounded operators, we have to show that $(1+A(D_N+i)^{-1})^{-1}$
  is bounded. We know in particular that $(D_N+i)^{-1}$ and therefore
  also $A(D_N+i)^{-1}$ is compact. Therefore $1+A(D_N+i)^{-1}$ is
  Fredholm of index $0$. Consequently, it is invertible if and only if its kernel
  is trivial. Now
  \begin{equation*}
    \begin{split}
      (1+A(D_N+i)^{-1}) f= 0 \quad&\iff  A (D_N+i)^{-1}f = - (D_N+i)
      (D_N+i)^{-1} f\\
      &\stackrel{g:=(D_N+i)^{-1}f}{\iff} (A+D_N)g = -i g.
    \end{split}
  \end{equation*}
  Since $D_N + A$ is selfadjoint, its spectrum does not contain $-i$ so
  that $g=0$ and therefore $\ker(1+ A(D_N+i)^{-1})=\{0\}$. Hence
  $(i+D_N+A)^{-1}$ and its summand $(i+D_V)^{-1}$ are
  $(2m+1)$-summable, too.

  For finite propagation speed we use the proof of
  \cite[Corollaire 1.11]{Hilsum(1989)}.
  There,   certain properties of the
  commutator $[D,h]$ with a Lipschitz function $h$ on $M$ are used.  Observe
  that the bundle homomorphism $A$ commutes with the multiplication
  operator $h$, therefore $[D_N,h]=[D_N+A,h]$ so that the proof for
  $D_N$ also applies to $D_N+A$. Since $D_V$ is a direct summand in
  $D_N+A$, finite propagation speed follows also for $D_V$. Exactly
  the same argument applies to $\overline{D_V}$.
\end{proof}

\begin{lemma}\label{lifttrclass}
  Let $R$ be a bounded trace class operator on $L^2(M,E)$ for
  some Lipschitz bundle $E$. Suppose 
  \begin{equation*}
\supp
  R(s)\subseteq U_\epsilon(\supp s):=\{x\in M\mid d(x,\supp(s))<\epsilon\} 
\end{equation*}
for every section $s\in L^2(M,E)$,
  and suppose the covering $\overline{M}\to M$ and the bundle $E$ are
  trivial over balls of radius $3\epsilon$. Then $R$ can be
  canonically lifted to
  a bounded operator $\overline{R}$ on $L^2(\overline{M},\overline{E})$ and
  $\overline{R}$ is of $\Gamma$-trace class with
  \begin{equation*}
    \tr_{\NeumannN\Gamma}(\overline{R})=\tr(R).
  \end{equation*}
\end{lemma}
\begin{proof}
  Decompose $M:=\disjointunion_{i=1}^n V_i$ with measurable subset $V_i$
  each of which has diameter less than $\epsilon$. Choose a lift
  $\overline{V_i}$ for each $V_i$. Then $\overline{M}
  =\bigcup_{i=1}^n\bigcup_{\gamma\in\Gamma}\gamma(\overline{V_i})$ and the
  union is disjoint up to sets of measure zero. Let $\phi_i^\gamma$ be the
  characteristic function of $\gamma(\overline{V_i})$.
  Every $\overline{s}\in
  L^2(M,E)$ is a sum $\sum \phi^\gamma_i \overline{s}$. By
  linearity we only have to define
  $\overline{R}(\phi_i^\gamma\overline{s})$
  $\forall i,\gamma$. We can identify the $2\epsilon$-neighborhood of
  $\gamma(\overline{V_i})$ with a corresponding neighborhood of $V_i$ in
  $M$, and since
  $R$ has only
  propagation $\epsilon$, in this way
  $\overline{R}(\phi_i^\gamma\overline{s})
  :=R(\phi_i^\gamma\overline{s})$ is well
  defined. Since $\abs{\overline{s}}_{L^2(\overline{M},\overline{E})}
  =\sum \abs{\phi_i^\gamma\overline{s}}_{L^2(\overline{M},\overline{E})}$
  and $R$ is bounded, this  makes sense also for the infinite sum
  $\sum \phi_i^\gamma\overline{s}$. In addition this 
show $\norm{\overline{R}}\le \norm{R}$.

  Let $\phi_i$ be the characteristic function of
  $V_i$. Multiplication with $\phi_i$ is a bounded operator on
  $L^2(M,E)$, therefore $R\phi_i$ is of trace class for each
  $i$. For fixed $i$, choose a fundamental domain of the covering
  which contains the $2\epsilon$-neighborhood of $\overline{V_i}$. This
  induces an obvious identification 
  \begin{equation*}
L^2(\overline{M},\overline{E})\cong
  L^2(M,E)\tensor l^2(\Gamma).
\end{equation*}
Moreover, under this identification
  the operator $\overline{R}_i = \sum_{\gamma\in\Gamma}
  \overline{R}\phi_i^\gamma$ becomes
  $R\phi\tensor \id_{l^2(\Gamma)}$. By standard properties of the
  $\Gamma$-trace (compare e.g.~\cite[Theorem 2.3(6)]{Schick(2001)}) 
  $\overline{R_i}$ is of $\Gamma$-trace class and
  $\tr_{\NeumannN\Gamma}(\overline{R_i})=\tr(R\phi_i)$
  (note that $\id\colon l^2(\Gamma)\to l^2(\Gamma)$ is of
  $\Gamma$-trace class with $\tr_{\NeumannN\Gamma}(\id)=1$). But $\overline{R}
=\sum_{i=1}^n\overline{R_i}$ and $R=\sum_{i=1}^n R\phi_i$. By
  linearity, $\overline{R}$ is of $\Gamma$-trace class with
  \begin{equation*}
   \hfill \tr_{\NeumannN\Gamma}(\overline{R})= \tr(R).\hfill\qed
  \end{equation*}
\renewcommand{\qed}{}
\end{proof}

Using the properties established in Theorem \ref{opprop} we can
essentially use Atiyah's proof to show:
\begin{theorem}
  In the situation described above
  \begin{equation*}
    \ind_{\NeumannN\Gamma}(\overline{D_V}) = \ind(D_V).
  \end{equation*}
  In particular, $\sign(M)=\sign^{(2)}(M)$.
\end{theorem}

We proceed with an outline of the proof. For details, we refer to Atiyah's
article \cite{Atiyah(1976)}. Assume throughout that $\dim M=4n$ is
divisible by four.

$D_V$ is an unbounded operator on $L^2\Omega^*(M,V)$,
$\overline{D_V}$ is its lift to
$L^2\Omega^*(\overline{M},\overline{V})$. The Hodge-$*$
operator induces a $\Z/2$-grading on $L^2\Omega^*$, and $D_V$ is an odd
operator with respect to this grading. What we are
really interested in is the graded index of $D_V$, i.e.~the index of
$D_V^+$ which maps the $+1$-eigenspace
of the grading operator
\begin{equation}\label{eq:def_of_tau}
\tau:=i^{p(p-1)+2n}*\qquad\text{(on $p$-forms)}
\end{equation}
to the
$-1$-eigenspace. Note that (using a fundamental domain)
$L^2\Omega^*(\overline{M},\overline{V})
\cong L^2(\Omega^*(M,V))\tensor l^2(\Gamma)$. The
problem is that kernel and cokernel of
$D_V$ and $\overline{D_V}$ can not be
related to each other using this product decomposition, because the corresponding
projection operators are highly nonlocal.

First step: construct a specific almost local parametrix for $D_V$ (the
same one is already used in \cite[Lemma 5]{Moscovici-Wu(1993)}). To
do this fix $\epsilon$ such that
the locally trivial covering $\overline{M}\to M$
 is trivial over balls of radius $3\epsilon$ (this is possible
since $M$ is compact). Choose a function $u\in C^\infty(\reals)$ such
that
\begin{enumerate}
\item $u$ is odd: $u(-x)=-u(x)$ $\forall x\in\reals$,
\item the function $v(x)= 1-xu(x)$ is rapidly decreasing,
\item the Fourier transforms of $u$ and $v$ are compactly supported
  with supports contained inside the interval $(-\epsilon,\epsilon)$.
\end{enumerate}
By Theorem \ref{opprop} \ref{selfadjoint} $D_V$ is selfadjoint. Using
functional calculus, we can construct $Q=u(D_V)$ and $R=v(D_V)$ and
conclude
\begin{equation*}
  D_VQ= 1 - R = QD_V.
\end{equation*}
Moreover, unit propagation speed
(see Theorem \ref{opprop} \ref{unitspeed}) implies that $Q$ and $R$ are
supported in an $\epsilon$-neighborhood of the diagonal,
i.e.~$\supp(Qf)\subseteq U_\epsilon(\supp(f))$ for any $f\in
L^2\Omega^*(M)$. By Lemma \ref{lifttrclass} we  can lift $Q$
and $R$ to operators $\overline{Q}$
and $\overline{R}$. Hence we lift the
whole equation to
\begin{equation*}
  \overline{D_V} \overline{Q} = 1-\overline{R}= \overline{Q}\overline{D_V}
\end{equation*}
(to check the that the domains coincide use that $\overline{D_V}$ is a
closed operator).

Second step: the parametrix property.
We required that $v$ is rapidly decreasing. This implies that $(i+x)^N
v(x)$ is bounded for every $N\in\naturals$ and therefore that $v(D_V) =
(i+D_V)^{-2m-1} \left((i+D_V)^{2m+1}v(D_V)\right)$
is of trace class, since by
Theorem \ref{opprop} \ref{summability} $(i+D_V)^{-1}$ is
$(\dim(M)+1)$-summable, therefore its $(\dim(M)+1)$st power is
$1$-summable, i.e.~of trace class.

Now remember that $D_V$ was anti-commuting with the grading operator
$\tau$ (i.e.~$\tau D_V=-D_V\tau$). Since
$u(x)$ is odd the same is true for $Q=u(D_V)$ by Lemma
\ref{grading} below. Since $v(x)=1-x u(x)$ is
even, $R=v(D_V)$ commutes with $\tau$. We therefore get a splitting
\begin{equation}\label{paramet}
 D_V^-Q^+= 1-R^+= Q^-D_V^+  ;\qquad D_V^+Q^- = 1- R^-= Q^+D_V^-
\end{equation}
where $R^{\pm}$ is the restriction of $R$ to the $\pm 1$-eigenspace of
$\tau$.
Since $\tau$ is a local operator, the operators $R^\pm$ are
$\epsilon$-local and their lifts are
$\overline{R}^{\pm}$. By Lemma \ref{lifttrclass}
$\overline{R}^{\pm}$ are of $\Gamma$-trace class and
\begin{equation*}
  \tr_{\NeumannN\Gamma}(\overline{R}^{\pm}) = \tr(R^{\pm}).
\end{equation*}

Step 3: Computing the index.
The main point now is that all the conditions are fulfilled to apply
Atiyah's principle of computing the index in terms of an arbitrary
parametrix. This is formalized in \cite[Proposition
2.6]{Schick(2001)}:  Let $H_0$
be the projection onto the kernel of $D_V^+$ and $H_1$ the projection
onto the cokernel of $D_V^+$ (which is the kernel of $D_V^-$ since $D_V^-$
is the adjoint of $D_V^+$). Define $T_0= (1-H_0) R^+ (1-H_0)$ and
$T_1= (1-H_1) R^- (1-H_1)$. Multiplication of \eqref{paramet} with
$H_0$ or $H_1$, respectively, yields $H_0=R^+H_0$ and
$H_1=H_1R^-$. This implies
\begin{align*}
  \tr(T_0) &=
  \tr(R^+)-\tr(H_0);\\
  \qquad\tr(T_1) &=\tr(R^-)-\tr(H_1).
\end{align*}
We want to show that $\ind(D_V^+)=\tr(H_0)-\tr(H_1)$ coincides
with $\tr(R^+)-\tr(R^-)$. To do this, it therefore suffices to
show that $\tr(T_0)=\tr(T_1)$. If $H$ is the projection onto
$\ker(D_V)$ then $D(1-H)R(1-H)=(1-H)R(1-H)D$ since all of these are
functions of $D$. Restriction to the positive subspace yields $
T_1D_V^+=D_V^+T_0$. Since $\ker(D_V^+)\subset
\ker(T_0)$ and $\ker(D_V^-)=\ker( (D_V^+)^*)\subseteq \ker (T_1)$,
$\tr(T_0)=\tr(T_1)$ is
the conclusion of \cite[Proposition 2.6]{Schick(2001)}
for the ordinary trace (where the group $\Gamma$ is trivial).

Exactly the same reasoning applies on the
universal covering $\overline{M}$ when computing the $\Gamma$-trace, to the
effect that
\begin{equation*}
  \ind_{\NeumannN\Gamma}(\overline{D_V}^+) = \tr_{\NeumannN\Gamma}(\overline{R}^+) -
  \tr_{\NeumannN\Gamma}(\overline{R}^-) =  \tr(R^+)-\tr(R^-) =\ind(D_V^+).
\end{equation*}

In the above proof we used:
\begin{lemma}\label{grading}
  Let $H$ be a $\integers/2$-graded Hilbert space with grading
  operator $\tau$. Let $D$ be a selfadjoint (not necessarily bounded)
  odd operator on $H$ (i.e.~$\tau D=-D\tau$). Let $f\colon \reals\to\reals$
  be a measurable function. If $f$ is odd or even then
  $f(D)$ is an odd or even operator, respectively.
\end{lemma}
\begin{proof}
  The grading operator is a unitary idempotent,
  i.e.~$\tau=\tau^*=\tau^{-1}$. Therefore $\tau^{-1}D\tau=
  -D$. Uniqueness of the spectral calculus implies $\tau^{-1}f(D)\tau=
  f(\tau^{-1}D\tau)$ for every function $f$. But $f$ even implies
  $f(-D)=f(D)$, and $f$ odd implies $f(-D)=-f(D)$ which concludes the proof.
\end{proof}

\begin{remark}\label{rem:Shmuels_proof}
  Shmuel Weinberger pointed out to us that one can also use a bordism
  argument to reduce the $L^2$-signature theorem for closed oriented topological
  manifolds to Atiyah's $L^2$-index theorem for closed oriented smooth manifolds.

  Indeed, every topological vector bundle $V$ over a topological
  manifold $M$ has a multiple which is topologically bordant to a
  smooth vector bundle over a smooth manifold (compare \cite[Theorem 1.6 and
  the following discussion]{Teleman(1984)}).

  It then remains to prove that the topological twisted signature is a
  bordism invariant. This is not clear from the classical proof of
  bordism invariance of the signature, which relies on the homological
  interpretation of the signature, and this is not available for
  twisted signature. However, Teleman \cite[Theorem 1.2]{Teleman(1984)} 
  proves the bordism invariance
  for the ordinary twisted signature, and we expect that a
  proof for the bordism invariance of twisted $L^2$-signature is
  possible along similar lines.
\end{remark}

When looking at manifolds with boundary, equality of signature and
$L^2$-signature fails as badly as possible. This follows from the fact
that essentially arbitrary intersection forms can be constructed, if the boundary
is non-empty. This is an easy consequence of Wall's non-simply
connected generalization of Milnor's plumbing
construction (compare \cite[Proof of Theorem 5.8]{Wall(1999)}). Since
we are not aware of a reference of this fact in
the literature, and since this is quite interesting a result, we prove
it here in reasonable detail.

\begin{proposition}\label{prop:arbitrary_signature_for_mf_with_boundary}
  Fix a dimension $2k\ge 6$ and a finitely presented group $\pi$. Let $X$ be a closed
  $(2k-1)$-dimensional manifold with fundamental group $\pi$ and with
  Morse decomposition without a $k$-handle. Let
  $V\iso (\integers\pi)^l$
  be a free finitely generated $\integers\pi$-module with (possible
  singular) $(-1)^k$-self dual map $\sigma\colon V\to
  V^*:=\Hom_{\integers\pi}(V,\integers\pi)$ of the form
  $\sigma=\psi+(-1)^k \psi^*$ (i.e.~$\sigma$ has a quadratic refinement).

  Then there is a compact manifold with boundary $(W; X, Y)$ of
  dimension $2k$ with boundary
  $\boundary W=X\disjointunion Y$ and with fundamental group $\pi$,
  such that the Morse chain complex
  $C_*(\tilde W)$ of the universal covering $\tilde W$ is isomorphic
  to $C_*(\tilde X)\oplus V$, where $V$ is considered as trivial chain
  complex concentrated in the middle dimension $k$, and with inverse
  Poincar{\'e} duality homomorphism
  \begin{equation*}
    C_{2k-*}(\tilde W) \to C_{2k-*}(\tilde W,\boundary\tilde
  W) \xrightarrow{PD^{-1}}
    C^*(\tilde W)
  \end{equation*}
  which in the middle dimension is exactly $\sigma$. Here $PD^{-1}$ is
  a chain homotopy inverse to the cup product with $[W,\boundary W]$.
\end{proposition}
\begin{proof}
  We use Wall's extension of Milnor's plumbing construction, as
  described in \cite[Proof of Theorem 5.8]{Wall(1999)}.

  More precisely, start with $X\times [0,1]$. Choose $l$ disjoint
  embedded discs $D^{2k-1}_i \subset X$. Let
  $i\colon S^{k-1}\times D^k\to D^{2k-1}$ be the standard
  embedding. By composition we obtain $r$ disjoint embeddings
  $f_i\colon S^{k-1}\times D^k_i\into X$. Choose
  lifts to the universal covering $\tilde X$. We now
  simultaneously deform the $f_i$ to new embeddings $f_i^1$ using regular
  homotopies $\eta_i$. The $\eta_i$ can be regarded as framed
  immersions of $S^{k-1}\times [0,1]$ to $X\times [0,1]$ (with
  boundary embedded). One can now count intersections and
  self-intersections as in \cite[(5.2)]{Wall(1999)} (taking the
  fundamental group into account using the chosen lifts). By
  \cite[p. 247]{Wall(1966)} the intersections and self-intersections
  can be chosen arbitrarily and independently.

  Now attach $k$-handles to $X\times [0,1]$ with attaching maps
  $f_i^1\times 1$. Let $W$ be the resulting manifold. Evidently,
  $\boundary W=X\disjointunion Y$, where $Y$ is obtained from $X$ by
  certain surgeries. Since the attaching maps are by construction
  homotopic to trivial embeddings, the statement about the cellular
  chain complex follows.

  It remains to adjust the intersection form. Choose the $\eta_i$ in
  such a way that the intersection of $\eta_i$ with $\eta_j$ is
  $\sigma(e_i)(e_j)$ where $\{e_i\}$ is the preferred bases of
  $(\integers\pi)^r$ and where we use the canonical isomorphism $V\iso
  {V^*}^*$. Moreover, choose $\eta_i$ such that the self-intersection
  of $\eta_i$ is $\psi(e_i)$. Then
  the intersection of $\eta_i$ with
  itself is $\psi(e_i)+ (-1)^k\psi(e_i)^*$, since our normal bundles
  are trivial.

  A canonical basis $\{S_i\}$ for the middle degree chain complex is given
  by the cores of the attached handles, completed to spheres using the
  images of the $\eta_i$ in $X\times [0,1]$ and the discs in
  $D^{2k-1}_i$ spanning the images of the $f_i$ (and with corners
  rounded). Then $S_i\cap S_j=\eta_i\cap \eta_j$, and the statement
  about the intersection form follows from the usual calculation of
  the Poincar{\'e} duality homomophism using intersection numbers.
\end{proof}

\begin{remark}
  Note that we could also prove a version of Proposition
  \ref{prop:arbitrary_signature_for_mf_with_boundary} for manifolds
  with middle dimensional handles in a Morse decomposition, with an
  additional summand in the middle degree chain complex.
\end{remark}

Observe that we use the usual translation of Poincar{\'e} duality to
homology, which, because of the use of intersection numbers is more
convenient to deal with in the case of smooth manifolds than the
cohomological version. Proposition
\ref{prop:arbitrary_signature_for_mf_with_boundary} implies, together
with Lemma \ref{lem:homology_intersection} the following corollary.

\begin{corollary}\label{corol:ord_and_L2_int_prescribed_for_mf_with_bound}
  If, in Proposition
  \ref{prop:arbitrary_signature_for_mf_with_boundary}, $X$ has a
  Morse decomposition without any $k$-cells  and $V\iso (\integers\pi)^l$,
  then the manifold $W$ has for an arbitrary $\integers\pi$-module $K$ 
  ``intersection form'' for homology twisted with $K$
  \begin{equation*}
    H_k(W;K) = K^l \xrightarrow{\id_K\tensor_{\integers\pi}\sigma}
    K^l\iso H^k(W;K). 
  \end{equation*}
   In particular, for
  $\pi$ the augmentation module $\epsilon\colon \integers\pi\to
  \reals$ we get the 
  ordinary intersection form
  \begin{equation*}
    H_k(W;\reals) = \reals^l
    \xrightarrow{\epsilon(\sigma)}
    \reals^l \iso H^k(W;\reals) = H_k(W;\reals)
  \end{equation*}
  where we use the canonical identification
  $H_k(W;\reals)=H^k(W;\reals)$ coming from cellular Hodge
  decomposition. Note that the (ordinary) signature is the signature
  of this self adjoint map (i.e.~the difference of the dimensions of
  positive and negative eigenspaces).

  Similarly, if $K=l^2(\pi)$ we get the ``$L^2$-intersection form'' 
  \begin{equation*}
    H^{(2)}_k(W) = (l^2(\pi))^l
    \xrightarrow{\sigma}
     (l^2(\pi))^l \iso H_{(2)}^k(W) = H^{(2)}_k(W)
  \end{equation*}
  where we use the canonical identification
  $H^{(2)}_k(W)=H_{(2)}^k(W)$ coming from cellular Hodge
  decomposition. Note that the $L^2$-signature is the $L^2$-signature
  of this self adjoint map (i.e.~the difference of $L^2$-dimensions of
  positive and negative spectral parts). Compare also
  \eqref{sign^(2) from L to R} and \eqref{sign from L to Z} and Section
  \ref{sec:comb-l2-sign}. 
\end{corollary}

  Note that, if $2k-1\ge 7$, for any finitely presented group $\pi$
  one can construct a closed manifold $X$ with fundamental group $\pi$
  and with a CW-structure without cells in dimension $k$.

\begin{theorem}
  Given any non-trivial finitely presented group $\pi$ and any
  dimension $4k\ge 8$, there is a manifold $W$ with boundary and with
  fundamental group $\pi$, such that 
  \begin{equation*}
    \sign^{(2)}(\tilde W)\ne \sign(W).
  \end{equation*}
\end{theorem}
\begin{proof}
  This follows immediately from Corollary
  \ref{corol:ord_and_L2_int_prescribed_for_mf_with_bound}, if we can
  produce appropriate (singular) intersection forms over
  $\integers\pi$. We use the fact that the signature and the
  $L^2$-signature can be computed in therms of the homology
  intersection form as well as the cohomological one, compare
  \ref{lem:homology_intersection}. 

  Any non-trivial group $\pi$ contains a non-trivial cyclic group
  $\Gamma$. Any finitely
  generated free $\integers\Gamma$ module with a given (possibly
  degenerate) intersection form can be induced up to a
  finitely generated free $\integers\pi$ module with induced intersection
  form, and the ordinary signature as well as the $L^2$-signature of
  the induced intersection form coincides with the original
  ones (compare also the proof of Remark \ref{rem: torsion-free Gamma is necessary for Poincare spaces}). Therefore, it suffices to treat the case $\pi$ cyclic.

  Using the canonical basis, we identify $(\integers\pi)^l$ with its
  dual. It suffices to consider the case $l=1$. In the case $\pi = \integers$ take $A$ to be the $(1,1)$-matrix $(1-z)$
for $z \in \integers$ a generator and let $\sigma\colon
  \integers\pi\to\integers\pi$ be given by multiplication with
  $A^*+A$. Then the augmentation
  $\epsilon\colon \integers\pi\to \reals$ gives $\epsilon(A +A^*) =
  0$ and yields zero as
  ordinary signature. The map $A+A^* \colon l^2(\integers) \to l^2(\integers)$ is a \emph{positive} weak
isomorphism and yields therefore the $L^2$-signature $1$ (the spectrum
  is contained in $[0,\infty)$, but there is no kernel). Notice that  
$A + A^*$ is not invertible over $\integers [\integers]$ so that we
  get no contradiction
to the conjecture that for torsion-free $\Gamma$ the maps
  $\sign^{(2)}$ and $\sign$ defined in 
  \eqref{sign^(2) from L to R} and \eqref{sign from L to Z} agree.

  If $\pi$ is a finite cyclic group of order $p>1$, we let $A=(1-z)$ and
  $\sigma\colon 
  \integers\pi\to\integers\pi$ again be given by multiplication with
  $(1-z) + (1-z^{-1})$ where $z$ is a generator of $\pi$. The
  augmentation yields the operator zero with ordinary signature
  $0$. On the other hand, on $l^2(\pi)=\complexs\pi$ the operator
  $A+A^*$ is non-negative with one-dimensional kernel (it diagonalizes
  with eigenvalues the values of $(1-z)+(1-z^{-1})$ at all $p$-th
  roots of unity). Therefore its signature (over $\complexs$) is
  $\dim_\complexs \complexs \pi -1 = p-1$. The $L^2$-signature is
  obtained by division by $\dim_\complexs \complexs \pi =p$ and
  therefore is $1-1/p \ne 0$.
\end{proof}


\section{$L^2$-index theorem for Poincar{\'e} spaces}
\label{sec:poincarecomp}

In this section we want to discuss special cases where the
$L^2$-signature theorem for closed Poincar{\'e} duality spaces is true. For finite
 fundamental groups, there are the counterexamples mentioned
in the introduction. For torsion-free fundamental  groups,
however, the $L^2$-signature theorem follows from the $C_{\smax}^*$-version of
the Baum-Connes conjecture or from the $L$-theory isomorphism
conjecture. 

Recall that there are symmetric $L$-groups $L^n_{\epsilon}(R)$ and quadratic $L$-groups
$L_n^{\epsilon}(R)$ for certain decorations $\epsilon =  p,h,s$ and $\langle - \infty \rangle$
and that there are symmetrization maps $L_n^{\epsilon}(R) \to L^n_{\epsilon}(R)$,
where in our context the ring with involution and unit $R$ is 
$\integers \Gamma$, $\rationals \Gamma$ or $\complexs\Gamma$
or the maximal group $C^*$-algebra $C_{\smax}^*\Gamma$. 
If one inverts $2$, then the decoration $\epsilon$ does not matter
and the symmetrization map is bijective. If we omit the decoration, we usually think of
$\epsilon = p$, i.e.\ the $L$-theory based on finitely generated projective modules. 
A reference for these definitions and facts is for instance
\cite[page 19, Section 1.10]{Ranicki(1981)}.  Note that for $C^*$-algebras $A$ there 
is a natural isomorphism between $L$-theory and topological $K$-theory
\cite[Theorem 1.6]{Rosenberg(1995)}
\begin{eqnarray} & L^n(A) \xrightarrow{\cong} K_n(A)& \label{L^n(A) = K_n(A)}
\end{eqnarray}
which will be used in the sequel without mentioning it. In dimension $n = 0$ it
sends the class of a non-degenerate sesquilinear form on 
a finitely generated projective module $P$
to the difference of the classes given by the positive part $P_+$ 
and by the negative part $P_-$.

The next definition is due to Wall \cite{Wall(1967)}:

\begin{definition} \label{def: Poincare pair} 
A \emph{$d$-dimensional Poincar{\'e} pair $(X,Y)$ over $\rationals$}
is a pair of finite  $CW$-complexes $(X,Y)$ such that  $X$ is connected, together with a
so called \emph{fundamental class} $[X,Y] \in H_d(X,Y;\rationals)$
such that for the universal covering and hence for any
$\Gamma$-covering $p\colon \overline{X} \to X$ the
Poincar{\'e} $\rationals \Gamma$-chain map induced
by the cap product with (a representative of) the fundamental class
$$- \cap [X,Y]\colon 
C^{d-*}(\overline{X},\overline{Y};\rationals) \to C_*(\overline{X};\rationals)$$
is a $\Q \Gamma$-chain homotopy equivalence.
If $Y = \emptyset$, we abbreviate $X = (X,\emptyset)$ and call it a 
\emph{$d$-dimensional Poincar{\'e} space}.
\end{definition}

Here $C_*(\overline{X};\rationals)$ is the cellular
$\rationals \Gamma$-chain complex and $C^{d-*}(\overline{X},\overline{Y};\rationals)$
is the dual $\rationals \Gamma$-chain complex
$\hom_{\rationals \Gamma}(C_{d-*}(\overline{X},\overline{Y};\rationals),\Q \Gamma)$.
Examples of Poincar{\'e} pairs are given by a compact connected
topological oriented manifold $X$ with boundary
$Y$ or merely by a connected closed oriented rational homology manifold.

\begin{theorem}\label{pcl2sig}
  Let $X$ be a $4n$-dimensional  Poincar{\'e} space over 
  $\rationals$. Let $\overline{X}\to X$ be a normal
  covering with \emph{torsion-free} covering group $\Gamma$. Assume
  that the (Baum-Connes) index map for the maximal group $C^*$-algebra
  \begin{equation*}
     \ind\colon K_0(B\Gamma) \to K_0(C^*_{\smax}\Gamma)
   \end{equation*}
   or the $L$-theory assembly map
 \begin{equation*}
     A\colon H_{4n}(B\Gamma;{\mathbb L}^{\langle - \infty \rangle}_{\bullet}) \to 
   L_0^{\langle - \infty \rangle} (\integers\Gamma)
   \end{equation*}
   is rationally surjective. Then
   \begin{equation*}
     \sign^{(2)}(\overline{X}) = \sign(X).
   \end{equation*}
\end{theorem}
\begin{proof}
  Since $X$ has no boundary, its symmetric signature 
  \begin{eqnarray}
  \sigma(\overline{X}) ~ \in ~ L^0(\rationals \Gamma)
  \label{symmetric signature}
  \end{eqnarray}  
  as an element
  in the symmetric projective $L$-group $L^0(\rationals\Gamma)$ is defined
  (for the definitions
  compare e.g.~\cite{Mishchenko(1970)}, \cite[Proposition 2.1]{Ranicki(1980bb)},
  \cite[page 26]{Ranicki(1981)}). 
  
  The $L^2$-signature $\sign^{(2)}(\overline{X})$ is the
  image of $\sigma(\overline{X})$ under the canonical map 
  \begin{eqnarray}
  \sign^{(2)} \colon L^0(\rationals \Gamma) & \to & \reals
  \label{sign^(2) from L to R}
  \end{eqnarray}
  which is the composition of change of rings homomorphism 
$L^0(\rationals \Gamma) \to L^0(\NeumannN\Gamma)$,
  the isomorphism  $L^0(\NeumannN\Gamma) = K_0(\NeumannN\Gamma)$ 
and the map induced by the standard trace
  $\tr_{\NeumannN\Gamma} \colon K_0(\NeumannN\Gamma) \to \reals$.
  The signature $\sign(X)$ is the
  image of $\sigma(\overline{X})$ under the canonical map 
  \begin{eqnarray}
  \sign \colon L^0(\rationals \Gamma) & \to & \integers
  \label{sign from L to Z}
  \end{eqnarray}
  which is the composition 
  $$L^0(\rationals \Gamma) \to L^0(\rationals) \to L^0(\complexs) =
  K_0(\complexs) = \integers.$$
  Hence it suffices to show that the maps $\sign^{(2)}$ and $\sign$ defined in 
  \eqref{sign^(2) from L to R} and \eqref{sign from L to Z} agree. 

   We begin with the case where the Baum-Connes index map is 
  assumed to be rationally surjective.
  By the Baum-Douglas description of $K$-homology, every element of
  $K_0(B\Gamma)$ is given by a map of a closed oriented smooth manifold $M\to B\Gamma$
  and an elliptic operator $D$ on $M$. Its index in
  $K_0(C^*_{\smax}\Gamma)$ is obtained by twisting $D$ with the pull
  back of the
  canonical $C^*_{\smax}\Gamma$-bundle on $B\Gamma$. The image of 
  this index element under the composition
  $$t^{(2)} \colon 
  K_0(C^*_{\smax}\Gamma) \to K_0(\NeumannN\Gamma) \xrightarrow{\tr_{\NeumannN\Gamma}} \reals$$
  can be read off directly as the $L^2$-index in
  the sense of Atiyah of the operator $\overline{D}$ lifted to the
  $\Gamma$-covering of $M$ which is the pull back of $E\Gamma$ via the
  map $M\to B\Gamma$. On the other hand, the image of this element under the composition
  $$t\colon K_0(C^*_{\smax}\Gamma) \to  K_0(C^*_{\smax}\{1\}) = K_0(\complexs) 
  \xrightarrow{\cong} \integers$$
  is just the index of $D$. (Here we need to deal with the maximal group $C^*$-algebra,
  because the reduced group $C^*$-algebra is not functorial under 
  group homomorphisms such as $\Gamma \to \{1\}$.)
  Atiyah's $L^2$-index theorem
  \cite[(1.1)]{Atiyah(1976)} now states
  that these two numbers coincide. Hence the two maps $t^{(2)}$ and $t$ above coincide since we assume
  that the index map $K_0(B\Gamma) \to K_0(C^*_{\smax}\Gamma) $ is rationally surjective.
  This implies that the maps $\sign^{(2)}$ and $\sign$ defined in
  \eqref{sign^(2) from L to R} and \eqref{sign from L to Z} above coincide since 
  $\sign^{(2)}$ and $\sign$, respectively, are given by the composition of 
  $t^{(2)}$ and $t$, respectively, with the map
  $$L^0(\rationals\Gamma) \to L^0(C_{\smax}^*\Gamma) \xrightarrow{\cong}
  K_0(C_{\smax}^*\Gamma).$$ 
  
  \medskip
  Now suppose that the $L$-theoretic assembly map is rationally surjective. 
  The symmetric signature defines for any $CW$-complex $Y$ a natural homomorphism
  $$\sigma \colon \Omega_*(Y) \to L^*(\integers\pi_1(Y)).$$
  The change of ring and decoration map
  $L_*^{\langle - \infty \rangle} (\integers \pi_1(Y)) \to L_*(\rationals \pi_1(Y))$ 
  and the symmetrization map $L_*(\rationals\pi_1(Y)) \to L^*(\rationals\pi_1(Y))$ are
  bijective after inverting $2$ \cite[pages 19, 104 and
  376]{Ranicki(1981)} and \cite[Proposition 8.2 or 3.3]{Ranicki(1980aa)}. 
  By the universal properties of assembly maps, 
  $\sigma  \otimes_{\integers} \rationals$ can be factorized as
  \begin{multline*}\sigma  \otimes_{\integers} \rationals \colon 
  \Omega_{*}(Y) \otimes_{\integers} \rationals 
  \to H_{*}(Y;{\mathbb L}_{\bullet}^{\langle - \infty \rangle}) \otimes_{\integers} \rationals  
  \xrightarrow{A \otimes_{\integers} \rationals}  
  L_*^{\langle - \infty \rangle}(\integers \pi_1(Y))\otimes_{\integers} \rationals  
  \\
  \xrightarrow{ \cong}  L^*(\rationals \pi_1(Y))\otimes_{\integers}
  \rationals .
  \end{multline*}
  where the first map is a transformation of homology theories with
  values in $\rationals$-vector spaces.
  The first map is surjective for $Y = \{Pt.\}$. Recall that every
  homology theory with values in rational vector spaces which vanishes
  in negative degrees is a direct
  sum of copies of shifted ordinary homology with rational coefficients
  (i.e.~the corresponding spectrum is a wedge of rational
  Eilenberg-Mac Lane spectra) (compare \cite{Dold(1968)}). It
  follows that the first map is surjective for all
$Y$. This could
  also be concluded using homological Chern characters. 
  The second map is surjective for $Y = B\Gamma$ by assumption and the
  third map is always bijective. Hence 
  $$\sigma \colon \Omega_*(B\Gamma) \to L^*(\rationals\Gamma)$$ 
  is rationally surjective. 
  This implies that rationally every element in $L^0(\rationals \Gamma)$ is a combination of
  elements of the form $\sigma(\overline{M})$ for $\Gamma$-coverings $\overline{M} \to M$
  with closed oriented smooth manifolds $M$ of dimension divisible by four as basis. 
  This follows also from the geometric interpretation of the
  assembly map in terms of the surgery sequence (see for instance
  \cite[Proposition 18.3]{Ranicki(1992)} for the topological category).
 For coverings $\overline{M}\to M$ as above we know already
  $\sign^{(2)}(\overline{M}) = \sign(M)$. Hence the maps  
  $\sign^{(2)}$ and $\sign$ defined in
  \eqref{sign^(2) from L to R} and \eqref{sign from L to Z} 
  above coincide. In \cite{Weinberger(1988b)}, a similar argument is used to prove
  homotopy invariance of $\rho$-invariants under the same assumptions
  we are making.
\end{proof}

   The ``max''-Baum-Connes conjecture used in Theorem \ref{pcl2sig} is
   true for $K$-amenable torsion-free
   groups for which the Baum-Connes conjecture is true, e.g.~torsion
   free amenable groups or torsion-free discrete subgroups of
   $SU(n,1)$ or $SO(n,1)$. For more information about 
   the Baum-Connes Conjecture see for instance
   \cite{Higson-Kasparov(2001)}, \cite{Mislin(2002)}, \cite{Valette(2002)}.

Examples of groups for which the $L$-theory isomorphism conjecture is
known are torsion-free poly-finite-or-cyclic groups
\cite{Farrell-Jones(1988b)}, fundamental groups of closed
non-positively curved manifolds \cite{Farrell-Jones(1993c)}, or knot
groups \cite{Aravinda-Farrell-Roushon(1997)}. 

\begin{remark} \label{rem: torsion-free Gamma is necessary for Poincare spaces}
We have seen in the proof of Theorem \ref{pcl2sig} that for a given finitely presented
group $\Gamma$ the $L^2$-index formula $\sign^{(2)}(\overline{X}) = \sign(X)$ holds for
each $\Gamma$-covering $\overline{X} \to X$ with a  $4n$-dimensional Poincar{\'e} space $X$ 
as base if the maps $\sign^{(2)}$ and $\sign$ defined in 
\eqref{sign^(2) from L to R} and \eqref{sign from L to Z} agree. It turns out that this
is an if and only if statement. Namely, rationally any element in $L^0(\rationals \Gamma)$ 
can be realized as $\sigma(\overline{X})$ by the following argument. Fix a closed manifold
$N$ of dimension $4n-1 \ge 7$ with $\pi_1(N) = \Gamma$ and 
$\eta \in L_0^s(\integers \Gamma)$. By Wall's realization theorem 
\cite[Theorem 5.8]{Wall(1999)} there is a normal map of degree one with underlying map 
$(f,\partial f) \colon (M,\partial M)
\to (N \times [0,1], N \times \{0,1\})$ such that $\partial f$ is a homotopy equivalence
and the associated surgery obstruction is $\eta$. The symmetrization map 
$L_0^s(\integers \Gamma) \to L^0_s(\integers \Gamma)$ sends the surgery obstruction to
the symmetric signature $\sigma(\overline{X})$ of the obvious $\Gamma$-covering of the
$4n$-dimensional Poincar{\'e} space $X$ which is
obtained by glueing $M$ and $N \times [0,1]$ together along their boundary with the
homotopy equivalence $\partial f$ \cite[Proposition 6.4]{Ranicki(1995b)}. Since the
composition
$$L_0^s(\integers \Gamma) \to L^ 0_s(\integers \Gamma) \to L^0(\integers \Gamma)  
\to L^0(\rationals \Gamma)$$
is bijective after inverting two, the claim follows.

It is not hard to check that the maps $\sign^{(2)}$ and $\sign$ defined in 
\eqref{sign^(2) from L to R} and \eqref{sign from L to Z} are different for $\Gamma$ 
a finite cyclic group of prime order (see for instance \cite[Example
22.28]{Ranicki(1992)}). Since for an inclusion $i \colon \Gamma \to \Gamma'$ of groups the
composition of the map $\sign^{(2)}$ for $\Gamma'$ with the induction homomorphisms
$i_* \colon L^0(\rationals \Gamma) \to L^0(\rationals \Gamma')$ is the
 map $\sign^{(2)}$ for $\Gamma$ and similar for $\sign$, the maps $\sign^{(2)}$ and
 $\sign$ for $\Gamma$ can only agree if and only if $\Gamma$ is torsion-free. 
In particular the conclusion in Theorem \ref{pcl2sig} that $\sign^{(2)}(\overline{X}) =
\sign(X)$ holds for $\Gamma$-coverings $\overline{X} \to X$ over 
$4n$-dimensional Poincar{\'e} spaces $X$  can only be true if $\Gamma$ is torsion-free.

\end{remark}

\begin{question}
To which extend does Theorem \ref{pcl2sig} hold for arbitrary
torsion-free groups?

Note that a negative answer would give rise to interesting elements in
the (quite mysterious and not well understood)
$K_0(C^*_{max}\Gamma)$ arising as (higher) signatures for closed
Poincar{\'e} duality spaces which, if the Baum-Connes conjecture for
$\Gamma$ is true, lie in the kernel of the map
$K_0(C^*_{max}\Gamma)\to K_0(C^*_r\Gamma)$. 
\end{question}


\section{Different definitions of $L^2$-signatures}
\label{sec:L2sig}

Throughout this section we consider a compact connected oriented $d = 4n$-dimensional Riemannian 
manifold $M$, possibly with boundary $\partial M$, 
together with a $\Gamma$-covering $\overline{M} \to M$.
We denote by $\overline{\partial M}$ the preimage of $\partial M$. More generally, we consider
a $d = 4n$-dimensional Poincar{\'e} pair
$(X,Y)$ over $\rationals$ with a $\Gamma$-covering $(\overline{X},\overline{X}) \to (X,Y)$. 
We will denote by $u\colon M \to B\Gamma$ and $u \colon X \to B\Gamma$ the classifying maps of
the $\Gamma$-coverings.

We present several different ways to define the
$L^2$-signature and show that they in fact coincide. 
One can use the $L^2$-index of the signature operator to define $\sign_{\san}^{(2)}(\overline{M})$
provided  $\partial M = \emptyset$. 
Using $K$-theory and  $L$-theory respectively  one can define $\sign_K^{(2)}(\overline{M})$ 
and $\sign_L^{(2)}(\overline{M})$  respectively
if $\partial M = \emptyset$. We will
define signature pairings on $L^2$-de Rham cohomology, and also on
combinatorial $L^2$-cohomology and take the von Neumann signature of
these. This will yield $\sign_{\sforms}^{(2)}(\overline{M},\overline{\partial M})$ 
and $\sign_{\schain}^{(2)}(\overline{X},\overline{Y})$.


\subsection{Analytic $L^2$-signatures}

\begin{definition}\label{sigop without boundary}
 
  Assume $\partial M=\emptyset$. The \emph{analytic $L^2$-signature}
  is the $L^2$-index (in the graded sense) of its signature operator,
  i.e.~if $\overline{D}=d+\delta$ is the signature operator on  $\overline{M}$
  and if
  $\overline{D}^\pm$ is its positive/negative
  part with respect to the signature
  splitting on $L^2\Omega^*(\overline{M})$ (i.e.~the restriction to the
  $\pm1$-eigenspace of $\tau=\pm *$ (compare \eqref{eq:def_of_tau})
  where $*$ is the
  Hodge-$*$-operator)
  then
  \begin{equation}\label{eq:def_of_L2_ind}
    \sign^{(2)}_{\san}(\overline{M}) := \ind_{\NeumannN\Gamma}(\overline{D}^+):=
    \dim_{\NeumannN\Gamma}(\ker\overline{D}^+ ) - \dim_{\NeumannN\Gamma}(\ker(\overline{D}^+)^*) .
  \end{equation}
  \end{definition}
  Note that $(\overline{D}^+)^*=\overline{D}^-$.
  This works not only for smooth Riemannian manifolds, but
  also for Lipschitz manifolds with Lipschitz Riemannian metrics and
  the corresponding Lipschitz signature operator.

  If $\partial M\ne \emptyset$ one still can use the signature
  operator. However, one has to supply it with the non-local
  Atiyah-Patodi-Singer boundary conditions. Moreover, to get the
  signature, one has to subtract a certain correction term
  (corresponding to ``extended $L^2$-solutions on the cylinder) from the
  index (compare
  \cite[(4.7)--(4.14)]{Atiyah-Patodi-Singer(1975a)}). To avoid this we
  directly \emph{define} the analytic index as the
  ``corrected cohomological'' expression of the index formula, namely,
   we put in the case $\partial M \not= \emptyset$
  \begin{eqnarray}
    \label{sigeq}
    \sign^{(2)}_{\san}(\overline{M},\overline{\partial M}) & := & \int_M L(M)
    -\eta^{(2)}(\overline{\partial M}) + \int_{\partial M} \Pi_L(\partial M).
  \label{sigop with boundary}
  \end{eqnarray}
  This coincides with
  the above definition if $\partial M=\emptyset$, and by
  \cite[Theorem 1.1]{Ramachandran(1993)} it also is the $L^2$-index of the
  signature operator (minus the standard correction term) if $\partial
  M\ne \emptyset$. 

\subsection{The $K$-theoretic $L^2$-signature}
\label{sec:K-theory_sig}

  Suppose $\partial M = \emptyset$. Form the flat twisted von Neumann algebra
  bundle $\mathcal{N}:=\NeumannN\Gamma\times_\Gamma\overline M$
with fiber the group von Neumann algebra $\NeumannN\Gamma$.
Given any elliptic differential operator $D\colon C^\infty(E)\to
  C^\infty(F)$ of order $d$ on $M$, one can twist this operator with
the bundle $\mathcal{N}$ to obtain an elliptic
  $C^*$-operator on $C^*$-vector bundles $\mathcal{E}$, $\mathcal{F}$.
An overview over
  this construction (for general $C^*$-bundles)
  can be found in
  \cite[Section 1]{Rosenberg(1983)}.

  One can define Sobolev spaces $H^s(\mathcal{E})$ of sections of
  $\mathcal{E}$, and similarly for $\mathcal{F}$. These are Hilbert
  $\NeumannN\Gamma$-modules,
  in particular, they have an inner product with values in
  $\NeumannN\Gamma$. The twisted operator then is a bounded operator
  \begin{equation*}
    D_{\mathcal{N}}\colon H^s(\mathcal{E})\to H^{s-d}(\mathcal{F}),
  \end{equation*}
  with a parametrix $Q\colon H^{s-d}(\mathcal{F})\to
  H^s(\mathcal{E})$.

  Then we define 
  \begin{equation*}
    \ind_{K_0(\NeumannN\Gamma)}(D_{\mathcal{N}}):=[\ker(D_{\mathcal{N}}+K)]
    -[\coker(D_{\mathcal{N}}+K)] \in K_0(\NeumannN\Gamma),
  \end{equation*}
  where we have to perturb by a $C^*$-compact operator $K$ to assure
  that kernel and cokernel are indeed finitely generated projective
  modules over $\NeumannN\Gamma$.

  The standard  trace $\tr_{\NeumannN\Gamma}$ defines (being a positive trace) a
  homomorphism
  \begin{equation*}
    \tr_{\NeumannN \Gamma} \colon K_0(\NeumannN\Gamma)\to \reals.
  \end{equation*}
  \begin{definition} \label{sign_K^{(2)}} If $\partial M = \emptyset$, 
  we define the \emph{$K$-theoretic $L^2$-index}
  \begin{equation*}
    \ind^{(2)}_{K}(D_{\mathcal{N}}):=\tr_{\NeumannN\Gamma}(
    \ind_{K_0(\NeumannN\Gamma)}(D_{\mathcal{N}})) \in\reals,
  \end{equation*}
  and the \emph{$K$-theoretic $L^2$-signature} as the corresponding index of
  the signature operator $D^+$:
  \begin{equation*}
    \sign_K^{(2)}(\overline M):= \ind_K^{(2)}(D^+_\mathcal{N}).
  \end{equation*}
\end{definition}

\begin{theorem}
  Suppose $\partial M = \emptyset$. For any elliptic differential operator $D$ on  $M$
  we have
  \begin{equation*}
    \ind_{\NeumannN\Gamma}(\overline{D}) = \ind^{(2)}_K(D_{\mathcal{N}}),
  \end{equation*}
  where $\overline D$ is the lift of $D$ to the $\Gamma$-covering,
  considered as unbounded operator on $L^2$-sections,
  and 
  \begin{equation*}
\ind_{\NeumannN\Gamma}(\overline{D}):=\dim_{\NeumannN\Gamma}(\ker(\overline
  D)) - \dim_{\NeumannN\Gamma}(\ker(\overline{D}^*))
\end{equation*}
is defined as in
  \eqref{eq:def_of_L2_ind} for the special case of the signature
  operator.

  In particular we get
  $$\sign_{\san}^{(2)}(\overline{M}) ~ = ~ \sign_K^{(2)}(\overline{M}).$$
\end{theorem}
  A proof for this well known result can be found in
  \cite{Schick(2002)}.

For the signature and the signature operator, the only operators we
are interested in here, we can actually rely on
a different set of results (already discussed at length in the
literature) which relate the higher signatures to surgery obstructions
in $L$-theory groups. This is discussed in Subsection
\ref{sec:L_group_sign}. 


\subsection{The de Rham $L^2$-signature}

Now we allow from the start that $\partial M\ne \emptyset$.

Let $V$ be a Hilbert space and let 
$s \colon V \times V \to \complexs$ be a sesquilinear pairing which is
bounded. For us, sesquilinear also means $s(v,w)=\overline{s(w,v)}$.
We can associate to it a selfadjoint bounded operator 
\begin{eqnarray}
&A \colon V  \to  V  &
\label{A associated to s}
\end{eqnarray} 
which is uniquely determined by the property
that $s(v_1,v_2) = \langle v_1,A(v_2)\rangle$ holds for all $v_1,v_2 \in V$.
From $A$ we obtain an orthogonal splitting $V = V_- \oplus V_- \oplus V_0$ 
of Hilbert spaces, where $V_+$ is the image of
$\chi_{(0,\infty)}(A)$, $V_0$ is the kernel of $A$ and  $V_-$ is the image of
$\chi_{(-\infty,0)}(A)$. The pairing $s$ is non-degenerate if and only
if $V_0$ is trivial. (One might want to require that $0$ is not
contained in the spectrum of $A$ as an ever stronger version of
non-degeneracy). 
If $V$ is a Hilbert module over the von Neumann algebra
$\NeumannN \Gamma$ and $s$ is $\Gamma$-invariant, then 
$A$ is $\Gamma$-equivariant and the splitting above
is a splitting of Hilbert $\NeumannN \Gamma$-modules. The \emph{$L^2$-signature} of
$s$ is in this case defined as 
\begin{eqnarray}
\sign^{(2)}(s) &  = & 
\dim_{\NeumannN \Gamma}(V_+) - \dim_{\NeumannN \Gamma}(V_-).
\label{von Neumann signature of a pairing}
\end{eqnarray}

  The cup-product of two $L^2$-forms is an $L^1$-form. If this
  product form is of the dimension of the manifold, we can
  integrate. In this way we get a pairing
  \begin{equation*}
  \innerprod{\cdot,\overline{\cdot}}\colon    L^2\Omega^p(\overline{M},\overline{\partial M})\times
    L^2\Omega^{4n-p}(\overline{M}, \overline{\partial M}) \to \complexs
  \end{equation*}
  which passes to $L^2$-cohomology as in the compact case. One
  should remark that this pairing factorizes through
  $\im(H^{p}_{(2)}(\overline M,\overline{\partial M})\to H^p_{(2)}(\overline M))$. 
  The restriction of the pairing to
  the middle dimension
  \begin{equation}
   s_{\sforms} \colon  H^{2n}_{(2)}(\overline{M},\overline{\partial M})\times H^{2n}_{(2)}
   (\overline{M},\overline{\partial M}) \to \complexs
  \label{formprod}  
  \end{equation}
  is a sesquilinear, bounded and $\Gamma$-invariant pairing.
  \begin{definition} \label{def: deRham signature}
  Define the \emph{de Rham $L^2$-signature} 
  \begin{equation*}
    \sign^{(2)}_{\sforms}(\overline{M},\overline{\partial M}) : =
    \sign^{(2)}(s)
  \end{equation*}
  to be the $L^2$-signature $\sign^{(2)}(s_{\sforms})$
  defined in \eqref{von Neumann signature of a pairing}
  for the pairing $s_{\sforms}$ introduced in 
  \eqref{formprod}.
\end{definition}

  Note that this does work for Lipschitz Riemannian manifolds
  as well as for smooth Riemannian manifolds.

  If $M$ is closed, the pairing is non-degenerate because to any $\omega\in
  L^2\Omega^{2n}(\overline{M})$ we can assign $*\omega\in
  L^2\Omega^{2n}(\overline{M})$ and
  $\int_{\overline{M}} \omega\wedge *\omega
  >0$ if $\omega\ne 0$. Moreover, we see that the splitting in this
  case is given by the $\pm 1$-eigenspaces of $*$: $H^+=\ker(*-1)$ and
  $H^-=\ker(*+1)$ (this makes sense if we identify the homology with
  the $L^2$-harmonic forms as can be done by Hodge theory). Moreover,
  the classical arguments apply to show that
  \begin{equation*}
    \ind_{\NeumannN\Gamma}(\overline{D}^+) = \dim_{\NeumannN\Gamma}(H^+)-\dim_{\NeumannN\Gamma}(H^-),
  \end{equation*}
  i.e.~all signatures $\sign_{\san}(\overline{M})$, $\sign_K(\overline{M})$ and
  $\sign_{\sforms}(\overline{M})$, defined so far, coincide. This also works for
  Lipschitz manifolds (compare \cite[Theorem 5.3]{Teleman(1983)}
  for the compact case).

The proof that $\sign_{\san}^{(2)}(\overline{M},\overline{\partial M})=
  \sign_{\sforms}^{(2)}
  (\overline{M}, \overline{\partial M})$ for manifolds with boundary
 (which (up to the usual error term) amounts to the fact that
  the index of the signature operator with APS-boundary conditions in
  fact gives the signature) is non-trivial even in the compact case, compare 
  \cite[(2.3)]{Atiyah-Patodi-Singer(1975a)} and
  the discussion after \cite[(4.5)]{Atiyah-Patodi-Singer(1975a)}. Moreover,
  this argument can not directly be used in the $L^2$-case, since it
  makes use e.g.~of a gap near zero in the spectrum of the signature
  operator on $\partial M$. To circumvent this requires considerable
  effort.

\begin{theorem}\label{the: L2indtheorem}
If  $M$ is a compact connected oriented $4n$-dimensional manifold with boundary $\partial M$ and
$\overline M \to M$ is $\Gamma$-covering as before, then
  \begin{equation}\label{sigeq3}
\sign_{\san}^{(2)}(\overline{M},\overline{\partial M})=
\sign_{\sforms}^{(2)}(\overline{M},\overline{\partial M}).
\end{equation}
\end{theorem}

First assume that the metric on $M$ has a product structure near the
   boundary. 
   The proof in the classical case in
   \cite{Atiyah-Patodi-Singer(1975a)} consists of two steps. In the first
   step they prove that the analytical index is the signature of the
   Poincar{\'e} duality pairing on the $L^2$-harmonic forms on
   $M_\infty$. Here $M_\infty$ is $M$ with an infinite cylinder $\partial
   M\times [0,\infty)$ attached to the boundary (with the product
   metric, which gives a smooth metric on all of $M_\infty$ because we
   started with a product metric near $\partial M$).

   We can similarly form $\overline M_\infty$ by attaching a cylinder
   to $\overline M$ (this is a $\Gamma$-covering of
   $M_\infty$). Let $\mathcal{H}_{(2)}^p(\overline M_\infty)$ be the
   $L^2$-harmonic $p$-forms on this manifold. Vaillant \cite[5.16]{Vaillant(1997)} 
   proves that the $L^2$-signature $\sign^{(2)}(s_{\sforms})$ of the intersection pairing 
   $$s_{\infty} \colon {\mathcal{H}}_{(2)}^{2n}(\overline M_\infty) \times  
   {\mathcal{H}}_{(2)}^{2n}(\overline M_\infty) \to \complexs.$$ 
   is $\sign^{(2)}_{\san}(\overline
   M)$.
   This is a non-trivial
   fact which we don't know a short and easy proof of. The $L^2$-version
   of the first step in the treatment in
   \cite{Atiyah-Patodi-Singer(1975a)} follows. Hence it remains to prove
   $$\sign_{\sforms}^{(2)}(\overline M) ~ = ~ \sign^{(2)}(s_{\infty}).$$
   We do this in the following sequence of lemmas.

   Remember first that we can define the $L^2$-homology of $\overline
   M$ as the reduced homology of the chain complex of
   $L^2$-differential forms on $\overline M$ (with no boundary
   conditions: compare \cite[Section 5]{Lueck-Schick(1998)} or
   \cite[Sections 1.4.2, 1.5]{Lueck(2002)} where a short
   account of different competing definitions is given).

Hence,
   restriction gives a map 
   \begin{equation*}
     r^p\colon{\mathcal H}^p(\overline M_\infty)\to H^p_{(2)}(\overline M).
   \end{equation*}
   We also have the natural map
   \begin{equation*}
     i^p\colon H^p_{(2)}(\overline M,\overline {\partial M})\to
     H^p_{(2)}(\overline M).
   \end{equation*}
   We will show that the closures of the image of $r^p$ and the image
   of $i^p$ coincide and
   that the pairings on ${\mathcal H}_{(2)}^{2n}(\overline M_\infty)$ and on
   $i^{2n}(H^{2n}_{(2)}(\overline M,\overline{ \partial M}))$ have the same $L^2$-signature.
Observe
that the pairing on $H_{(2)}^{2n}(\overline M,\overline{\partial M})$ 
is well defined by a standard integration by parts argument, and
the same argument shows that it descends to
$\im(i^{2n}\colon H_{(2)}^{2n}(\overline M,\overline{\partial M})\to
H_{(2)}^{2n}(\overline M))$. 

 We first
   prove:
\begin{lemma}\label{correct range}
     The image of $r^p$ lies in the closure of the image of $i^p$. 
 \end{lemma}
\begin{proof} Let
   $$q^p \colon H^p_{(2)}(\overline{M}) \to H^p_{(2)}(\overline{\partial M})$$
   be the map given by restriction. To prove the statement, because of the long weakly exact
   sequence  for the $L^2$-cohomology of the pair $(\overline M,\overline{\partial M})$
   we only have to check that $q^p \circ r^p$ vanishes.
   If $\omega\in {\mathcal H}(\overline M_\infty)$ then by definition
   $\omega$ is $L^2$-integrable. Because of elliptic regularity, it
   lies in $H^\infty:=\bigcap_{s\ge 0} H^s$, i.e.~all derivatives
   are in $L^2$.
   In particular, using the continuous restriction homomorphism to
   codimension $1$ submanifolds $H^s(\overline{M}_\infty)\to
   L^2(\boundary \overline{M}\times \{t\})$ ($s>1/2$), for $t \in [0,\infty)$ the
   pull back map indeed gives $L^2$-forms on 
  $\overline{\partial M} \times \{t\} = \overline{\partial M}$, i.e.
   \begin{equation*}
     q[t]^p\colon{\mathcal H}^p(\overline M_\infty) \to L^2\Omega^p(\overline{\partial M}).
   \end{equation*}
   Notice that $q[0]^p = q^p \circ r^p$. 
   The maps  $q[t]^p$ are continuous, and all the manifolds 
   $\overline{\partial M}\times [r,\infty)$ are isometric. Given a form $\omega\in
   {\mathcal H}^p(\overline M_\infty)$, the sequence of its
   restrictions $\omega_t$ to $\overline{\partial M}\times [t,\infty)$
   tends to zero in all Sobolev norms (where we use the isometry just
   described to compare the different $\omega_t$). Therefore the sequence 
   $q[t]^p(\omega)$ in $L^2\Omega^p(\overline{\partial M})$ tends to
zero as $t\to \infty$.

   Now all forms
   $q[t]^p(\omega)$ represent the same element in the reduced
   $L^2$-homology of $\overline{\partial M}$. This is true since, on
   the cylinder, we can write $\omega=\omega_1(u)+\omega_2(u)\wedge du$ (if
   $u$ is the cylinder variable), with $\omega_{1,2}$ $L^2$-functions
   on $[0,\infty]$ with values in $L^2\Omega^*(\overline{\partial M})$. Observe that $\omega$ is
   closed. Therefore
   \begin{equation*}
     0 =d\omega = d\omega_1(u) \pm \frac{\partial\omega_1(u)}{\partial
     u}\wedge du + d\omega_2(u)\wedge du
 \end{equation*}
 Since the summands with and without $du$ are linearly independent,
 from this we get
 \begin{equation*}
   \pm \frac{\partial\omega_1(u)}{\partial
     u}  = (d\omega_2(u)).
 \end{equation*}
Integrating this equation with respect to $u$ we get
\begin{equation*}
  \omega_1(t)-\omega_1(0) = \pm d(\int_0^t\omega_2(u) \;du).
\end{equation*}
But $\omega_1(t)$ is the pullback of $\omega$ to the submanifold
$\overline{\boundary M}\times \{t\}$, and we conclude
   \begin{equation}\label{eq:1}
     q[t]^p(\omega)-q[0]^p(\omega)= \pm d(\int_0^t\omega_2(u) \;du).
   \end{equation}

   We consider $\omega_{1,2}$ to be $L^2$-functions on $[0,\infty)$
   with values in the Hilbert space $L^2\Omega^*(\overline{\boundary
   M})$. To those, we can apply the Cauchy-Schwarz inequality: the
   inner product of $\omega_2(u)$ and the constant function with value
   $1$ satisfies:
   \begin{multline}
\abs{\innerprod{\omega_2(u),1}_{L^2([0,t]; L^2\Omega^*(\overline{\tilde
   M}))}}^2 = \abs{\int_0^t \omega_2(u)\;du}^2\\
\le
   \int_0^t 1^2 \;du\cdot \int_0^t\abs{\omega_2(u)}^2\;d u \le
   t\int_0^t \abs{\omega(u)}^2\;du, 
 \end{multline}
Therefore the difference on the left hand side of Equation
   (\ref{eq:1}) is the differential of an
   $L^2$-form. Because $q[t]^p(\omega)\xrightarrow{t\to\infty} 0$ in
   $L^2$, this proves the lemma.
\end{proof}

From here one, we cannot continue exactly as in the classical case,
because forms
representing zero are not exactly boundaries, and homology sequences
are only weakly exact. Instead, we use von Neumann dimensions and
suitable subspaces with codimensions tending to zero.

First we address surjectivity of the restriction map 
\begin{equation*}
r^p \colon 
{\mathcal{H}}_{(2)}^p(\overline M_\infty)\to \overline{\im (i^p)} = 
\ker q^p \subseteq H^p_{(2)}(\overline M).
\end{equation*}

Consider the differential $d\colon\Omega^{p-1}_{(2)}(\overline{\partial M})\to
\Omega^p_{(2)}(\overline{\partial M})$. This map is unbounded and left Fredholm by
elliptic regularity (compare e.g.~\cite[Lemma 3.3]{Lueck-Schick(1998)})
and hence the image of the spectral projection $\chi_{(0,\gamma)}(\delta
d)$ has von Neumann dimension which tends to zero for
$\gamma\to0$. For given $\epsilon>0$ choose $\gamma > 0$ such that the
image of $\chi_{(0,\gamma)}(\delta d)$ has dimension not bigger than
$\epsilon$. Put 
$$E_{\epsilon}^p~:= ~ \im(d \circ \chi_{(\gamma,\infty)}(\delta d)) \subseteq
\Omega^p_{(2)}(\overline{\partial M}).$$
 Since $d \circ \chi_{(-\infty,0]}(\delta d)$ is zero, 
$E_{\epsilon}^p$ has codimension $\le\epsilon$ in $\overline{\im
(d)}$. Moreover, $E_{\epsilon}^p$ is closed since the restriction of
$\delta d$ to the relevant subspace fulfills $\delta d\ge\gamma$ and
hence is invertible.

If, using the well established Hodge decomposition (compare
e.g.~\cite[Theorem 5.10]{Schick(1996)url})
\begin{multline} \label{Hodge decomposition old}
  L^2\Omega^{2n-1}(\overline{M})  = \\
  \overline{(\im(\overline{d}^{2n-2}))}
 \oplus
 \overline{(\im(\overline{\delta}^{2n}|_{\{\omega ;
   \omega|_{\overline{\partial M}=0}\}}))} \oplus
 \ker(\overline{\Delta}_{2n-1}|_{\{\omega;
   (*\omega)|_{\overline{\partial
  M}}=0=(\delta\omega)|_{\overline{\partial M}}\}}),
 \end{multline}
we identify $H^p_{(2)}(\overline M)$ with the space of harmonic
forms which fulfill absolute boundary conditions, pulling back to the
boundary gives a well
defined bounded map $H^p_{(2)}(\overline M)\to
\Omega^p_{(2)}(\overline{\partial M})$. Let 
$$K_{\epsilon}^p \subseteq H^p_{(2)}(\overline M)$$
be the inverse image of $E_{\epsilon}^p$ under this map. It is a closed subspace of 
$H^p_{(2)}(\overline M)$ which actually is contained in $\ker(q^p)$, 
the inverse image of $\overline{\im(d(\overline{\partial M}))}$, and has codimension $\le
\epsilon$ in $\ker(q^p)$.

\begin{lemma}\label{surjectivity}
  $K_{\epsilon}^p$ is contained in the image of 
$r^p\colon {\mathcal{H}}_{(2)}^p(\overline M_\infty)\to \overline{\im(i^p)}$.
\end{lemma}
\begin{proof}
  Let $\omega$ be a harmonic form representing an element in
  $K_{\epsilon}^p$. Then we have to find a harmonic form 
  $h\in{\mathcal{H}}_{(2)}^p(\overline M_\infty)$ whose restriction to
  $\overline{M}$
  represents the cohomology
  class of $\omega$. By assumption, $q^p\omega=d\alpha$ for suitable
  $\alpha\in \Omega^{p-1}_{(2)}(\overline{\partial M})$ in the domain
  of $d$. Note that $d\alpha$ itself is smooth by elliptic regularity
  since $\omega$ is harmonic. Choose a
  smooth function $\psi\colon[0,\infty) \to\reals$ with  $\psi(t)=1$ in a
  neighborhood of $0$ and with $\psi(t)=0$ for $t>1/2$. Define
  $\tilde\alpha=\alpha\cdot \psi(t)\in
  \Omega^{p-1}_{(2)}(\overline{\partial M}\times 
  [0,\infty))$. Note that $\tilde{\alpha}$ is an $L^2$-form in the
  domain of $d$, which is smooth in a neighborhood of the
  boundary. For such forms, all usual integraion by parts formulas
  hold, a fact we are using frequently in the sequel and which follows
  e.g.~from the methods of \cite{Gaffney(1954)}, or is explained in
  more detail in \cite{Schick(1996)url}.

  Let $Q\colon\overline{\partial M}\times\{0\}\into \overline{\partial
  M}\times [0,\infty)$ 
  be the inclusion. Then $Q^{p-1}\tilde \alpha = \alpha$ and $Q^p
  d\tilde\alpha = d\alpha$.
  Define the $L^2$-form
  $\tilde\omega$ on $\overline M_\infty$ to coincide with $\omega$
  on $\overline M$, and with $d\tilde\alpha$ on 
  $\overline{\partial M}\times[0,\infty)$. We claim that $\tilde\omega\in\ker(d)$,
  i.e.~that $\tilde\omega$ is orthogonal to $\delta\phi$ for all
  smooth $\phi$ with compact support. This is checked by integration
  by parts (on $\overline M$ and $\overline{\partial M}\times
  [0,\infty)$ separately): since $\omega$ is closed
  \begin{equation*}
    \innerprod{\tilde{\omega}|_{\overline{M}}, \delta\phi|_{\overline{M}}}_{L^2(\overline M)} =
    - \int_{\overline{\partial M}} d\alpha\wedge q[0]^{4n-1-p}(*\overline{\phi}),
  \end{equation*}
  on the other hand
  \begin{equation*}
    \innerprod{\tilde{\omega}|_{\overline{\partial M}  \times [0,\infty)},
\delta\phi|_{\overline{\partial M} 
\times [0,\infty)}}_{L^2(\overline{\partial M}\times [0,\infty))} 
   = - \int_{\overline{\partial M}^-} d\alpha\wedge
    q[0]^{4n-1-p}(*\overline{\phi}). 
  \end{equation*}
  Because of opposite inward
  directions the orientation of $\overline{\partial M}$ in the first and second
  integral are different. Changing the orientation changes 
the sign of the integral of a differential form.
   This implies the vanishing of
  $\innerprod{\tilde\omega,\delta\phi}_{L^2(\overline M_\infty)}$,
  which is just the sum of the two terms above.

  By Hodge decomposition, we therefore can
  write $\tilde \omega=h+x$ where
  $h\in{\mathcal{H}}_{(2)}^p(\overline{M}_\infty)$ and $x$
  lies in the closure of the image of $d$. If we apply $r^p$ to this
  equation, we see
  that the forms $\omega$ and $r^p(h)$ represent the same
  $L^2$-cohomology class in $H^p_{(2)}(\overline{M})$,  which finishes the proof. 
\end{proof}

\begin{corollary} \label{cor: dense image of r^p}
  The map $r^p\colon {\mathcal{H}}_{(2)}^p(\overline M_\infty)\to \overline{\im(i^p)}$ 
has dense image.
\end{corollary}
\begin{proof}
  The map surjects onto subspaces of arbitrary small codimension.
\end{proof}

Now we have to compare the intersection forms. Again we can not do
this directly, but have to restrict our attention to subspaces with
small codimension. Observe that $q[0]^p$ defines a map from
${\mathcal{H}}_{(2)}^p(\overline M_\infty)$ to $\overline {\im
  d(\overline{\partial M})}$. Let
$${\mathcal{H}}^p_\epsilon \subseteq {\mathcal{H}}_{(2)}^p(\overline M_\infty)$$ 
be the inverse image of $E_\epsilon$ under
this map. On the space of harmonic forms, the pull back map is bounded
in the $L^2$-norm, therefore ${\mathcal{H}}^p_\epsilon$ is closed. The
codimension of ${\mathcal H}^p_\epsilon$ in ${\mathcal{H}}_{(2)}^p(\overline
M_\infty)$ is not
bigger than $\epsilon$.

\begin{lemma} \label{lem: fors are almost equal}
  Let $\omega,\eta\in {\mathcal{H}}_\epsilon^{2n}$, with
  $q[0]^{2n}\omega=d\alpha$ and $q[0]^{2n}\eta=d\beta$. Define $\tilde
  \alpha$ and $\tilde \beta$ as in the proof of Lemma
  \ref{surjectivity}. Assume, without loss of generality, that
  $\overline M$ has a collar of length $1$ which is
  isometric to a product. Define $\tilde \alpha'$ and $\tilde
  \beta'$ as above, but with support on this collar of $\overline
  M$ (i.e.~replacing the ``outward'' $\boundary M\times [0,\infty)$ by
  the ``inward'' collar). Then $v:=r^{2n}(\omega)-d\tilde \alpha'$ and
  $w:=r^{2n}(\eta)-d\tilde\beta'$ pull back to zero on
  $\overline{\partial M}$ and represent the same homology classes as
  $r^{2n}(\omega)$ and $r^{2n}(\eta)$, respectively. Moreover,
  \begin{equation}\label{formequal}
    \int_{\overline M} v\wedge w = \int_{\overline M_\infty} \omega\wedge\eta.
  \end{equation}

\end{lemma}
\begin{proof}
  We only have to prove Equation \ref{formequal}. Integration by parts
  shows that
  \begin{equation}\label{eq:alpha_dalpha}
    \int_{\overline M}v\wedge w=\int_{\overline M}\omega\wedge\eta +
    \int_{\overline{\partial M}} \alpha\wedge d\beta,
  \end{equation}
since the additional terms $\int_{\overline{\boundary
  M}}d\alpha\wedge\alpha$ and $\int_{\overline{\boundary
  M}}d\beta\wedge\beta$ vanish as $2d\alpha\wedge\alpha=d\alpha\wedge
  \alpha+\alpha\wedge d\alpha = d(\alpha\wedge\alpha)=0$ as $\alpha$
  is of odd degree.

  We therefore have to show that
  \begin{equation*}
    \int_{\overline{\partial M}\times [0,\infty)}\omega\wedge\eta =
    \int_{\overline{\partial M}} \alpha\wedge d\beta.
  \end{equation*}

  Write $\omega-d\tilde\alpha = h_1+x$ and
  $\eta-d\tilde\beta=h_2+y$, where we restrict to
  $\overline{\partial M}\times [0,\infty)$ and use the Hodge
  decomposition for closed forms with vanishing pullback to the
  boundary. This implies that the harmonic forms $h_1$ and $h_2$ also
  fullfill $q^{2n}(h_1)=0=q^{2n}(h_2)$, and
  $x,y\in\overline{d(\im(d_b))}$, where $d_b$ stands for the
  differential $d$, but with domain only the smooth compactly
  supported forms whose pull back to the boundary is zero. Integration
  by parts shows that
  \begin{equation*}
    \int_{\overline{\partial M}\times [0,\infty)} (h_1+x)\wedge(h_2+y)
    = \int_{\overline{\partial M}\times [0,\infty)} h_1\wedge h_2.
  \end{equation*}
  We can write $h_1=a(t)+b(t)\wedge dt$ and $h_2=c(t)+b(t)\wedge dt$,
  and because of the product 
  structure the fact that $h_1$ is harmonic implies that the form $a$
  is harmonic and the form $b$ (or equivalently $b\wedge dt$) is
  harmonic. But $0=q^{2n} h_1 = a(0)$, and a harmonic form which
  vanishes identically at the boundary is zero, therefore $a=0$. In
  the same way, $c=0$. This implies $h_1\wedge h_2=0$ since $dt\wedge
  dt=0$. Consequently
  \begin{equation*}
    0 = \int_{\overline{\partial M}\times [0,\infty)}
    (\omega-d\tilde\alpha)\wedge(\eta-d\tilde \beta) =
    \int_{\overline{\partial M}\times[0,\infty)} \omega\wedge\eta -
    \int_{\overline{\partial M}} \alpha\wedge d\beta,
  \end{equation*}
  where the last equation follows from integration by parts (see
  \cite{Gaffney(1954)} as in \eqref{eq:alpha_dalpha}.
  This finishes the proof of Lemma \ref{lem: fors are almost equal}.
\end{proof}

Let 
$$L_{\epsilon}^{2n} \subseteq \overline{\im(i^p)}$$
be the closure of the image of ${\mathcal{H}}^{2n}_\epsilon$ under 
$r^p\colon {\mathcal{H}}_{(2)}^p(\overline M_\infty)\to \overline{\im(i^p)}$. The codimension of
$L_{\epsilon}^{2n} \subseteq \overline{\im(i^p)}$ is $\le \epsilon$ because of 
Corollary \ref{cor: dense image of r^p}, since the codimension of
${\mathcal H}^p_\epsilon$ in ${\mathcal{H}}_{(2)}^p(\overline
M_\infty)$ is $\le \epsilon$. 
The intersection form 
$$s_{\schain} \colon H^{2n}_{(2)}(\overline{M},\overline{\partial M}) \times 
H^{2n}_{(2)}\overline{M},\overline{\partial M})  \to \complexs$$
descends to a pairing on $\overline{\im(i^p)}$ which can be restricted to a paring
$$s \colon L_{\epsilon}^{2n}  \times L_{\epsilon}^{2n} \to \complexs.$$
Since the codimension of $L_{\epsilon}^{2n} \subseteq \overline{\im(i^p)}$ is $\le \epsilon$,
we get 
$$|\sign^{(2)}(s_{\schain}) - \sign^{(2)}(s)| ~ \le ~ \epsilon.$$
Lemma \ref{lem: fors are almost equal} implies that
the intersection form 
$$s_{\infty} \colon {\mathcal{H}}_{(2)}(\overline{M}_\infty) \times  
{\mathcal{H}}_{(2)}(\overline{M}_\infty) \to \complexs$$
restricts to a pairing on ${\mathcal{H}}^{2n}_{\epsilon}$ which descents to
the pairing $s \colon L_{\epsilon}^{2n}  \times L_{\epsilon}^{2n} \to \complexs$ above.
Since the codimension of
${\mathcal H}^p_\epsilon$ in ${\mathcal{H}}_{(2)}^p(\overline
M_\infty)$ is $\le \epsilon$ we get
$$|\sign^{(2)}(s_{\infty}) - \sign^{(2)}(s)| ~ \le ~ \epsilon.$$
We conclude  
$$|\sign^{(2)}(s_{\infty})  -  \sign^{(2)}(s_{\schain})| \le 2\epsilon.$$
Since $\epsilon > 0$ was arbitrary,
we get 
$$\sign^{(2)}(s_{\infty})  =  \sign^{(2)}(s_{\schain}).$$
This finishes the proof of
Theorem \ref{the: L2indtheorem} in the case, where the 
Riemannian metric is a product metric near $\partial M$.

  The argument also shows that $r^p\colon
  {\mathcal{H}_{(2)}^p}(\overline M_\infty)\to H^p_{(2)}(\overline M)$ is
  injective. This is the
  case because the intersection pairing
  is non-degenerate on ${\mathcal{H}}_{(2)}^p(\overline M_\infty)$ (if $0\ne
  h\in {\mathcal{H}}_{(2)}^p(\overline M_\infty)$
  then $h$ is not perpendicular to $*h$ where $*$ is the Hodge
  operator), and because on subspaces of
  arbitrarily small codimension this passes to the image of $r^p$.

  The general version of Theorem \ref{the: L2indtheorem}  (without
  product metric near the boundary) now follows by observing that
  $H^*_{(2)}(\overline M,\overline{\partial M})$ 
  is unchanged if we deform the
  metric on $M$ to a product metric, and that the intersection form
  also does only
  depend on the homology. We can deform the metric in such a way that
  the restriction to the boundary is unchanged (but of course the
  second fundamental form changes). If one does this, in
  $\sign_{\san}^{(2)}(\overline M)$ only the local terms
  $\int_{ M} L( M)$ and $\int_{\partial M}
  \Pi_L(\partial M)$ are changed. Exactly the same changes appear in
  the Atiyah-Patodi-Singer index formula for the ordinary signature on
  the compact manifold
  $M$. We know that the classical index formula also for manifolds without product
  metric near the boundary computes the signature, which does not
  depend on the metric. Therefore the overall changes are zero, and
  the same is true for $\sign_{\san}^{(2)}(\overline M)$. Since we
  just argued that $\sign_{\sforms}^{(2)}(\overline M)$ does not depend
  on the metric on $M$,  Theorem \ref{the: L2indtheorem} follows. \qed



\subsection{The combinatorial $L^2$-signature}
\label{sec:comb-l2-sign}

  Now we want to give a combinatorial construction of the
  pairing in \eqref{formprod}. Assume
  therefore that instead of a compact connected oriented Riemannian
  manifold $M$ we have a $4n$-dimensional Poincar{\'e} pair $(X,Y)$ over $\rationals$. 
  Recall that the Poincar{\'e} structure is given by
  a fundamental class $[X,Y] \in H_{4n}(X,Y;\rationals)$ with the
  following property.   Let the fundamental chain  
  $[X,Y] \in C_{4n}(  X, Y;\rationals)$, denoted in the same way
  as the fundamental class, be a closed chain representing the fundamental class.
  Let $(\overline{X},\overline{Y}) \to (X,Y)$ be
  a regular covering.
  Lift this closed chain $[X,Y]$ to the covering $\overline{X}$.
  The lift will be a closed bounded chain (without compact
  support) $[\overline{X},\overline{Y}]
  \in L^\infty C_{4n}(\overline{X},\overline{Y})=
   l^\infty(\Gamma)  \tensor_{\complexs \Gamma}
  C_{4n}(\overline{X},\overline{Y};\complexs)$. There is a duality
  pairing $L^1C^m(\overline{X},\overline{Y};\complexs)\times
  L^\infty C_m(\overline{X},\overline{Y};\complexs) \to \complexs$.
  We call the pairing against $[\overline{X},\overline{Y}]$
  ``integration over $\overline{X}$''. Now the cup product of one
  $L^2$-cochains with the complex conjugate of a second one on
  $\overline{X}$ gives an $L^1$-cochain which, if the
  dimensions are right, can be integrated over $\overline{X}$.
  This passes to reduced $L^2$-(co)homology 
  \begin{definition}\label{combdef} Denote the induced sesquilinear
  $\Gamma$-invariant bounded pairing of Hilbert $\NeumannN\Gamma$-modules (in the middle dimension
  $2n$) by
  \begin{equation}\label{eq:pairing}
   s_{\schain} \colon H_{(2)}^ {2n}(\overline{X},\overline{Y}) \times 
   H_{(2)}^{2n}(\overline{X},\overline{Y}) \to \complexs.
  \end{equation}
  Define the \emph{combinatorial $L^2$-signature}
  $\sign^{(2)}_{\schain}(\overline{X},\overline{Y})$ to be the
  associated $L^2$-signature $\sign_{\NeumannN\Gamma}(s_{\schain})$ of
  $s_{\schain}$ 
  as in \eqref{von Neumann signature of a
  pairing}.
  \end{definition}

  To show that this definition makes sense, recall that the definition of the
  cup-product involves a cellular approximation to the diagonal
  embedding $X\to X\times X$, which we can lift to an equivariant
  cellular map $\overline{X}\to \overline{X}\times \overline{X}$. This way, there
  is a global bound $K$ such that the image of each cell in $\overline{X}$
  under the diagonal approximation meets only $K$ cells of 
  $\overline{X}\times \overline{X}$. Remember that the cochain representing the
  cup-product of $a$ and $b$ maps a cell $\sigma$ to a certain linear
  combination of $a(\sigma_1)\cdot b(\sigma_2)$ n (given locally
  by the diagonal approximation), where
  $\sigma_1\times\sigma_2$ runs through all cells in the image of
  $\sigma$ under the cellular approximation to the diagonal. This
  implies in a standard way that this cup-product map is continuous
  from the product of the $L^2$-cochain spaces to the $L^1$-cochains.

  The result of the pairing between a cochain $\sum_{\sigma\text{
      $p$-cell}}\lambda_\sigma\sigma$  and a chain of the form $\sum_{\sigma\text{ 
      $p$-cell}}\mu_\sigma\sigma$ is the number
  $\sum_{\sigma\text{ $p$-cell}} \lambda_\sigma\overline{\mu_\sigma}$. This
  is a continuous pairing between $L^1$-cochains and
  $L^\infty$-chains. 

  Taken together, we get a pairing on $L^2$-cochains with values in
  the complex numbers. If we restrict in one factor to cochains with compact
  support, this is the classical pairing. In particular, $\int_{\overline{X}}
  a\cup b=0$ if $a=\delta(a')$ and $a'$ has compact support and
  $\delta(b)=0$, since this is true (in the classical situation) if 
  $a$ has compact support and $b$ is completely arbitrary. We want to
  check the corresponding statement if $a$ is in the closure of the
  image of $\delta$ in the space of $L^2$-cochains, and $b$ is an
  $L^2$-cochain with $\delta(b)=0$. Now $a=\lim_{n\to\infty}
  \delta(a_n)$, where we can assume that all $a_n$ have compact
  support, because $\delta$ is continuous and the cochains 
  with compact support are dense in the space of $L^2$-cochain. But
  then continuity implies the claim that our pairing vanishes on (the
  closure of the space of)
  coboundaries and therefore passes to reduced $L^2$-cohomology. The
  usual proofs apply to show that the cup product (and the pairing)
  does not depend on the particular way we constructed it (e.g.~the
  particular cellular approximation to the diagonal
  embedding).

Note that the
  construction is homological in nature and therefore depends only on
  the oriented homotopy type of the pair $(X,Y)$. In particular it is
  independent of the CW-structure and the choice of the closed cycle representing
  the fundamental class.

  An alternative description of Definition \ref{combdef} can be given
  using the sequence
  \begin{equation}\label{eq:chain}
    C^{4n-*}_{(2)}(\overline{X},\overline{Y})\xrightarrow{-\cap 
 [\overline{X},\overline{Y}]}
    C^{(2)}_*(\overline{X})\to C^{(2)}_*(\overline{X},\overline{Y}) .
  \end{equation}
  Note that this is obtained by tensoring the corresponding
  $\complexs\Gamma$-chain map over $\complexs\Gamma$ with $l^2(\Gamma)$. 
  It induces a selfadjoint bounded $\Gamma$-equivariant operator  
  \begin{equation}
    A\colon H^{2n}_{(2)}(\overline{X},\overline{Y}) \to H_{2n}^{(2)}(\overline{X},\overline{Y})
  \xrightarrow[\iso]{g}  H^{2n}_{(2)}(\overline{X},\overline{Y})
  \label{Poincare version of A} 
  \end{equation}
  using the canonical identification $H_{2n}^{(2)}(\overline{X},\overline{Y})
  =  H^{2n}_{(2)}(\overline{X},\overline{Y})$ which comes from the cellular Hodge 
  decomposition. Actually, putting any positive inner product on
  $H^{(2)}_{2n}(\overline{X},\overline{Y})$ will give rise to an
  identification with its dual space
  $H^{2n}_{(2)}(\overline{X},\overline{Y})$, and the fact that the
  Poincar{\'e} duality homomorphism is self dual implies that after the
  identification the homomorphism is self adjoint (with respect to the
  used inner product), as can be seen by going through the definitions.

\begin{lemma}\label{lem:homology_intersection}
  The homological Poincar{\'e} duality homomorphism
  \begin{equation*}
    B\colon H^{(2)}_{2n}\overline{X}) \xrightarrow[\iso]{PD^{-1}}
  H^{2n}_{(2)}(\overline{X},\overline{Y})
  \xrightarrow{i^*}H^{2n}_{(2)}(\overline{X})
  \xrightarrow[\iso]{g^{-1}} H^{(2)}_{2n}(\overline{X})
  \end{equation*}
  has the same $L^2$-signature as $A$ of \eqref{Poincare version of
  A}, where $PD^{-1}$ is the defined to be the inverse of the
  isomorphism induced by cup product with the fundamental class (which
  we abbreviate with $PD$ in this lemma). The
  corresponding remark holds for the ordinary signature of $(X,Y)$.

  For the calculation of $L^2$-signatures and ordinary signatures, the
  Poincar{\'e} duality chain map can be replaced by any chain homotopic
  map, and moreover, it can be ``conjugated'' with a chain homotopy
  equivalence and its adjoint.
\end{lemma}
\begin{proof}
  Given any self-adjoint Hilbert $\NeumannN\Gamma$-module morphism
  $a\colon V\to W$ and a (not necessarily unitary) Hilbert
  $\NeumannN\Gamma$-isomorphism $f\colon V\to W$, we have
  \begin{equation*}
    \sign^{(2)}(a) = \sign^{(2)}(faf^*).
  \end{equation*}
  This follows from the fact that the isomorphism $f^*$ intertwines
  $a$ and $faf^*$, i.e.
  \begin{equation*}
    \innerprod{(faf^*)x,x} =\innerprod{a (f^*x),(f^*x)}\qquad\forall
    x\in W,
  \end{equation*}
  i.e.~$f^*$ maps the positive or negative spectral part,
  respectively, of $faf^*$ to the corresponding part of $a$, and being
  an $\NeumannN\Gamma$-isomorophism, it preserves the
  $\NeumannN\Gamma$-dimension.

  In our case, 
  \begin{multline*}
\sign^{(2)}(B) = \sign^{(2)}( PD^* \circ B\circ PD) = \sign^{(2)}(PD^*
g^{-1} i^*)\\
=\sign^{(2)}((g^{-1}\circ i_*\circ PD)^*) =\sign^{(2)}(A),
\end{multline*}
since first $i^*\colon H^{2n}_{(2)}(\overline X,\overline Y)\to
H^{2n}_{(2)}(\overline X)$ is dual to $i_*\colon
H_{2n}^{(2)}(\overline X)\to H_{2n}^{(2)}(\overline X,\overline Y)$,
and therefore $i^*\circ g^{-1}$ is adjoint to $g^{-1}\circ i_*$ by the
usual relations between dual and adjoint on Hilbert spaces, and
secondly $A^*=A$.

The statement about the chain homotopy invariance follow trivially
from the fact that $L^2$-signature and signature depend on the
homological Poincar{\'e} duality map only, which is not affected by
passing to a chain homotopic map, and ``conjugation'' with a chain
homotopy equivalence and its adjoint corresponds to ``conjugation'' by
an isomorphism and its adjoint. We have just checked that this does not
change the $L^2$-signature.

The identical argument applies to the ordinary signature (which can be
considered as the $L^2$-signature for the trivial one-sheeted covering).
\end{proof}

  The standard relations between cup- and cap-product and
  ``integration'' of homology against cohomology classes imply
  \begin{proposition}\label{PD_map}
  The operator in \eqref{Poincare version of A} is the operator associated in 
  \eqref{A associated to s} to the pairing appearing in Definition \eqref{combdef}.
  In particular we get
  \begin{equation}\label{eq:chain_def_of_L2_sig}
    \sign^{(2)}_{\schain}(\overline{X},\overline{Y}) =
    \dim_{\NeumannN\Gamma}(\chi_{(0,\infty)}(A)) - 
\dim_{\NeumannN\Gamma}(\chi_{(-\infty,0)}(A)).
  \end{equation}
  \end{proposition}

  Suppose $(X,Y)$ happens to be an oriented cocompact smooth manifold
  with boundary, and the
  CW-structure is given by a smooth triangulation. Then by the $L^2$-de
  Rham isomorphism of Dodziuk \cite[Theorem 1]{Dodziuk(1977)}
  and its version
  for manifolds with boundary (\cite[Corollary 1.7]{Schick(1996)url} or
  \cite{Goldstein_et_al(1988)}),
  $L^2$-simplicial and $L^2$-de Rham cohomology are isomorphic.

For reasons of completeness, we will prove that the pairings which give
rise to $\sign^{(2)}_{\sforms}$ and $\sign^{(2)}_{\schain}$ are compatible 
with respect to this isomorphism. It would perhaps be more
satisfactory to prove that the isomorphism is compatible with the
products. However, we don't want to discuss the $L^1$-version of the
Hodge-de Rham theorem (and note that the product of two $L^2$-forms is 
an $L^1$-form), so we use this shortcut. The advantage is that we can
give a ``local'' proof of the weaker result, which holds on the chain
level. Note that, in contrast, there is no good way to describe a good cup
product on the level of cochains of a simplicial complex which is
at the same time graded commutative and associative, as is the
case for the wedge product of differential forms. Similar and related
work has e.g.~been done in \cite[Section 7]{Miller(1998)}, and his
methods could be used as well. Another version would use an
intermediate simplicial $L^2$-de Rham complex as in the treatment
in \cite{Dupont(1978)} of multiplicativity of the ordinary de Rham
isomorphism. Actually, this method is used in
\cite{Goldstein_et_al(1988)} to prove the de Rham theorem for
$L^2$-cohomology (as well as $L^p$-cohomology), but without taking
care of the multiplicative structure. We believe that the combination
of \cite{Dupont(1978)} and
\cite{Goldstein_et_al(1988)} proves that the $L^2$-de Rham
isomorphism preserves the multiplicative structure.

We choose to give a direct argument, using some calculations of
\cite{Ranicki-Sullivan(1976)}.

To start with, we recall a possible definition of the cup product on
the cochain level of a simplicial complex (using the Alexander-Whitney 
approximation).

So, assume $X$ is a simplicial complex. Choose an orientation of $X$,
i.e.~an orientation of each simplex of $X$. Next, we choose a local
ordering of the chain
complex, i.e.~a total ordering of the vertices of every simplex with
the compatibility condition that,
if a simplex $\sigma$ is the face of a simplex $\tau$, then the
restriction of the ordering on the simplices of $\tau$ should give the 
ordering on $\sigma$. Customarily, such a local ordering is obtained by 
globally ordering all the vertices of the simplicial complex, but that 
is by no means necessary for the following cup product construction,
and for us it will later be much more convenient to use local
orderings.

Observe that we do not require that the ordering is compatible with
the orientation (later on, we will use different local orderings, but
the same orientation).

If $e_0,\dots,e_n$ are the ordered vertices of a simplex $\sigma$, then 
$\langle e_0,\dots,e_n\rangle:=\epsilon(e_0,\dots,e_n)\sigma$ is a
chain, where $\epsilon(e_0,\dots,e_n)=1$ if $(e_0,\dots,e_n)$
represents the orientation of $\sigma$, and
$\epsilon(e_0,\dots,e_n)=-1$, otherwise.

Following the conventions in \cite{Rinow(1975)}, the cup product 
of a $p$-cochain $a$ and a $q$-cochain $b$ is
defined by
\begin{equation}\label{eq:def_of_cup_product}
  a\cup b(\langle e_0,\dots,e_n\rangle)=a(\langle e_0,\dots,e_p\rangle))\cdot
  b(\langle e_{p},\cdots,e_n\rangle ).
\end{equation}

Note in particular that, if $a$ is the elementary cochain corresponding to
$\langle e_0,\dots,e_p\rangle$ (i.e.~maps this simplex to one, and all 
other simplices to zero), and $b$ is the elementary
cochain of $\langle e_p,\dots, e_n\rangle$, then $a\cup b$ is
the elementary
cochain of $\langle e_0,\dots,e_n\rangle$.

The de Rham map $\int$ maps a (sufficiently smooth) $p$-form $\omega$
to a $p$-cochain of the
simplicial cochain complex of a smooth triangulation of the
manifold. The value of $\int(\omega)$ on a $p$-simplex $\sigma$ simply is
the integral of $\omega$ over $\sigma$. This is a chain map.

An inverse map $W$ from the cochain complex to differential forms
(going back to Whitney) is given by mapping an elementary $p$-cochain
$\sigma$ with vertices $(e_0,\dots,e_p)$ to the ``barycentre form''
\begin{equation*}
  W(\sigma):= p! \sum_{i=0}^p (-1)^i x_i\;dx_0\wedge\cdots\wedge
  \widehat{dx_i}\wedge \cdots \wedge dx_p,
\end{equation*}
where the hat means, as usual, that the corresponding entry is omitted,
and the $x_i$ are defined to be the barycentric coordinates  with
non-zero values in the stars of the vertices $e_i$. The form
$W(\sigma)$ is non-zero only on the open star of $\sigma$.

Dodziuk \cite{Dodziuk(1977)}, compare also
\cite{Goldstein_et_al(1988)}, proves that $W$  indeed induces 
an isomorphism on reduced $L^2$-cohomology. The inverse is essentially
induced by $\int$. In particular, it is easily established that $\int\circ
W=\id$. However, since $\int$ is (below the top degree) not
defined for all $L^2$-forms, one has to be somewhat careful here. This is
the main reason why we don't prove that the de Rham isomorphism is
multiplicative for $L^p$-cohomology (where the product of an $L^p$ and
an $L^q$-form is an $L^r$-form with $1/r=1/p+1/q$).

\begin{lemma}\label{lem:pairing_and_integration}
  If $c_1$ and $c_2$ are elements of the simplicial $L^2$-cochain
  complex such that the 
  degrees add up to $4n$, then on the
  $4n$-dimensional manifold $X$
  \begin{equation*}
    s_{\schain}(c_1,c_2) = \sum_{\sigma\text{oriented $4n$ simplex of $X$}}
    (c_1\cup \overline{c_2})(\sigma) = \int_X W(c_1\cup \overline{c_2}),
  \end{equation*}
  where the sum is over all $4n$-simplices with orientation induced
  from $X$.
\end{lemma}
\begin{proof}
  The first equality is the definition of the pairing. For the second
  one observe that
  \begin{multline*}
      \int_X W(c_1\cup \overline{c_2}) = \sum_{\sigma\text{ $4n$-simplex of $X$}}
      \int_\sigma W(c_1\cup \overline{c_2})\\
      = \sum_{\sigma\text{ $4n$-simplex of $X$}} (\int\circ W)(c_1\cup
      \overline{c_2})(\sigma) 
      = \sum_{\sigma\text{oriented $4n$-simplex of $X$}} (c_1\cup
      \overline{c_2})(\sigma). 
    \end{multline*}
    We used the fact that $\int\circ W$ is identically the identity map.
\end{proof}

Since we already know that $W$ induces an isomorphism, from this it suffices to
check for the compatibility of the two pairings that 
\begin{equation}\label{eq:product_compatibility}
 \innerprod{W(c_1),W(c_2)}:=  \int_X W(c_1)\wedge
 \overline{W(c_2)}  = \int_X W(c_1\cup \overline{c_2}) 
\end{equation}
for $c_1$ and $c_2$ cochains as in Lemma
\ref{lem:pairing_and_integration}, since the right hand side equals
$s_{\schain}(c_1,c_2)$. 

We are only interested in the result on cohomology. Therefore, we can define
the cup product on the cochain level appropriately. Recall, as already
observed above, that many choices are possible. Our description
depends e.g.~on the chosen local ordering.

First, we assume that our triangulation is the barycentric subdivision
of some other triangulation (if it is not yet, pass to the barycentric
subdivision). There, a canonical local ordering is defined: a simplex $\sigma$
of the barycentric subdivision is by definition a chain $s_0\subset
s_1\subset\dots\subseteq s_k$ of simplices of the original triangulation
with vertices $s_0$,\ldots, $s_k$; and the ordering on the latter is
given by inclusion, or, equivalently, by ordering according to the dimension.

In the latter description, on our $4n$-dimensional simplicial complex
$X$ we define a collection of local orderings parameterized by the
symmetric group $\Sigma_{4n+1}$ of permutations of $\{0,\dots,4n\}$
with vertices $s_i$, $s_j$ of the simplex $\sigma$ above satisfying
$s_i<_\tau s_j$ under the ordering induced by $\tau\in
\Sigma_{4n+1}$ if and only if $\tau(s_i)<\tau(s_j)$.We denote the cup
product induced by this local ordering by $\cup_\tau$.

The cup product to be used for Equation
\eqref{eq:product_compatibility} is then the average of all the
$\cup_\tau$:
\begin{equation*}
  c_1\cup c_2:= \frac{1}{(4n+1)!} \sum_{\tau\in \Sigma_{4n+1}}
  c_1\cup_\tau c_2. 
\end{equation*}

We now prove Equation~\eqref{eq:product_compatibility} with this
definition of the cup product.
\begin{proof}
  Let $v_1(c_1,c_2):= \int_X W(c_1)\wedge \overline{W(c_2)}$ and
  $v_2(c_1,c_2):=\int_X W(c_1 \cup \overline{c_2})$ for simplicial $L^2$-cochain
  $c_1,c_2$.

  Then $v_1$ and $v_2$ are sesquilinear and jointly continuous. For the
  latter we use the fact that $W$ is a continuous map from
  $L^2$-cochain to $L^2$-forms as well as from $L^1$-cochain to
  $L^1$-forms (this follows from its ``local'' character). Moreover,
  the wedge as well as our cup product are continuous from $L^2$ to
  $L^1$ by an appropriate application of the H{\"o}lder inequality
  (again, the ``local'' definition of the cup product is used here).

  The span of the elementary cochains given by the (oriented) simplices of the
  triangulation (defined after Equation \eqref{eq:def_of_cup_product})
  is dense in the space of all $L^2$-cochains. Consequently, it
  suffices to prove that $v_1(c_1,c_2)=v_2(c_1,c_2)$ if $c_1$ and
  $c_2$ are two cochains corresponding to oriented simplices
  $\sigma_1= (e_0,\dots,e_p)=\langle e_0,\dots,e_p\rangle$ or 
  $\sigma_2=(f_0,\dots,f_q)=\langle f_0,\dots,f_q\rangle$, respectively.

  Let us first consider the case that $\sigma_1$ and
  $\sigma_2$ have no vertex in common. Then the cup product of $c_1$
  and $c_2$ is
  zero. At the same time, the supports of $W(c_1)$ and $W(c_2)$ (being
  the open stars of the simplices $\sigma_1$ and $\sigma_2$) have
  empty intersection. In this case therefore
  $v_1(c_1,c_1)=0=v_2(c_1,c_2)$.

  Secondly, assume $\sigma_1$ and $\sigma_2$ have $2$ or more vertices
  in common. Then $W(c_1)\cup \overline{W(c_2)}=0$ since each summand
  contains the square of a one-form $dx_j$ for some barycentre
  function $x_j$. Similarly, $c_1\cup \overline{c_2}(\sigma)=0$ for each
  \emph{non-degenerate} simplex and in particular for each (non-degenerate)
  $4n$-simplex, so again $v_1(c_1,c_2)=0=v_2(c_1,c_2)$.

  Finally, for the interesting case, assume $f_0=e_p$ is the only
  vertex which both simplices have in common ($f_0=e_p$ is no real loss of
  generality, we could replace an oriented simplex by the negative of
  a simplex with the wrong orientation and the whole argument would go
  through). Evidently, only the case
  $p+q=4n$ is of interest, in
  which case $(e_0,\dots,e_p=f_0,f_1,\dots f_q)$ spans a
  $4n$-simplex. Let $\sigma$ be the oriented simplex with these
  vertices and with orientation induced from $X$. Observe that
  $\sigma$ is spanned by $\sigma_1$ and $\sigma_2$, but the
  orientation it gets that way differs from its orientation by
  $\epsilon(e_0,\dots,f_q)=:*(\sigma_1,\sigma_2)$. The latter notation
  is used in \cite[p.~23]{Ranicki-Sullivan(1976)}.

The support of
  $W(\sigma_1)\cup \overline{W(\sigma_2)}$ is the
  interior of this $4n$-simplex.
  Therefore, its integral over any other $4n$-simplex is zero.

  Moreover, $c_1\cup c_2$ vanishes on all $4n$-simplices apart from
  $\sigma$ (as follows immediately from the formula for the cup
  product), and hence $W(c_1\cup c_2) = (c_1\cup c_2)(\sigma)$. It
  remains to compute this number. Our definition of the cup product
  involves one summand for each of the $(4n+1)!$ permutations of the
  simplices of $\sigma$. The contribution of such a permutation can
  only be nontrivial, when
  the first $p+1$ simplices $(e_0,\dots,e_p)$ are mapped to themselves and
  the last $q+1$-simplices $(f_0,\dots,f_q)$ are also mapped to themselves, in particular, $e_p=f_0$
  has to be fixed by such a permutation. Observe that we obtain exactly
  $p!\cdot q!$ permutations with non-trivial contribution.

  \begin{multline}\label{eq:cup_calc}
    c_1\cup \overline{c_2}(\langle e_0,\dots,e_p=f_0,\dots, f_q\rangle) =\\
 \frac{1}{(p+q+1)!}
    \sum_{\pi\in \Sigma_p,\psi\in\Sigma_q} c_1(\langle
    e_{\pi(0)},\dots, e_{\pi(p-1)}, e_q\rangle) \overline{c_2(\langle f_0,
    f_{\psi(1)},\dots f_{\psi(q)}\rangle)}. 
  \end{multline} 
  The definition of the chain 
  \begin{equation*}
\langle
  e_{\pi(0)},\dots,e_{\pi(p-1)},e_p\rangle
\end{equation*}
differs from the simplex
    $(e_{\pi(0)},\dots,e_{\pi(p-1)},e_p)$ by a sign which makes up for
    the (possible) change of orientation compared to the oriented
    simplex spanned by $e_0,\dots,e_p$. This implies that the value of the expression in
    \eqref{eq:cup_calc} does not depend on the particular
    permutation. For our cup product, we therefore get
    \begin{multline*}
      c_1\cup \overline{c_2}(\langle e_0,\dots,e_p=f_0,\dots, f_q\rangle) =\\
    \frac{p!\cdot q!}{(p+q+1)!} c_1(\langle e_0,\dots,e_p\rangle )
    c_2(\langle f_0,\dots, f_q\rangle) =  \frac{p!\cdot q!}{(p+q+1)!}
  \end{multline*}
  by the definition of $c_1$ and $c_2$.
  Finally, observe that
  \begin{equation*}
    c_1\cup \overline{c_2}(\sigma) =*(\sigma_1,\sigma_2) c_1\cup \overline{c_2} (\langle
    e_0,\dots, f_q\rangle) =*(\sigma_1,\sigma_2) \frac{p!\cdot q!}{(1+p+q)!}.
  \end{equation*}

  It remains to calculate $\int_\sigma W(c_1)\wedge
  \overline{W(c_2)}$. This is carried out in
    \cite[Appendix]{Ranicki-Sullivan(1976)} and we obtain indeed
    \begin{equation*}
      \int_\sigma W(c_1)\wedge \overline{W(c_2)} = *(\sigma_1,\sigma_2)
      \frac{p!\cdot q!}{(p+q+1)!}.
    \end{equation*}

    This finishes the proof of the claim.
\end{proof}

In particular, it follows that:
\begin{proposition}
   Assume $M$ is a compact oriented smooth manifold with boundary $\partial M$. Then
   \begin{equation*}
       \sign^{(2)}_{\schain}(\overline{M},\overline{\partial M}) = 
      \sign^{(2)}_{\sforms}(\overline{M},\overline{\partial M}).
   \end{equation*}
\end{proposition}


\subsection{The $L$-theoretic $L^2$-signature}
\label{sec:L_group_sign}

\begin{definition}\label{L_theory_def}
  Consider a Poincar{\'e} space $X$ of dimension $d = 4n$ over $\rationals$.
Let $\overline{X} \to X$ be a regular $\Gamma$-covering.  
We have already mentioned its symmetric signature
$\sigma(\overline{X}) \in L^0(\integers\Gamma)$ in \eqref{symmetric signature}.
Define its \emph{$L$-theoretic $L^2$-signature}
$$\sign_L^{(2)}(\overline{X}) ~ \in ~ \reals$$
as the image of $\sigma(X)$ under the map $\sign^{(2)} \colon L^0(\rationals \Gamma) \to \reals$
introduced in \eqref{sign^(2) from L to R}. 
\end{definition}

\begin{lemma}\label{lem:L_theory_sig_and_L2chain_complex_sig}
  In the situation of Definition \ref{L_theory_def}, we have
  $$\sign_L^{(2)}(\overline{X}) ~ = ~ \sign_{\schain}^{(2)}(\overline X). $$
\end{lemma}
\begin{proof} 
Let $\mathcal{U}(\Gamma)$ be the algebra of operators affiliated to $\NeumannN \Gamma$.
Algebraically  $\mathcal{U}(\Gamma)$ is the Ore localization 
of $\NeumannN \Gamma$ and has the property
that it is a von Neumann regular ring, i.e. 
any finitely generated submodule of a finitely generated projective
$\mathcal{U}\Gamma$ module is a direct summand \cite[Theorem 8.22]{Lueck(2002)}.
There is a commutative square
$$\begin{CD}
L^0(\NeumannN \Gamma) @>>> K_0(\NeumannN \Gamma)
\\
@VV\cong V @VV \cong V
\\
L^0(\mathcal{U} \Gamma) @>>> K_0(\mathcal{U} \Gamma)
\end{CD}
$$
where the vertical maps are change of rings maps and isomorphisms
\cite[Theorem 9.31]{Lueck(2002)}. Since $\mathcal{U} \Gamma$ is von
Neumann regular, the $\mathcal{U}\Gamma$-chain complex 
\begin{equation*} 
\cdots \xrightarrow{0} H_*(C_*(\overline{X};\rationals) \otimes_{\rationals \Gamma}
\mathcal{U}\Gamma)\xrightarrow{0} \cdots
\end{equation*}
given by the homology and the trivial differentials
consists of finitely generated projective $\mathcal{U} \Gamma$-chain modules and
there is a $\mathcal{U}\Gamma$-chain 
homotopy equivalence 
$$i_* \colon H_*(C_*(\overline{X};\rationals) 
\otimes_{\rationals \Gamma} \mathcal{U}\Gamma)
\to C_*(\overline{X};\rationals) \otimes_{\rationals \Gamma} \mathcal{U}\Gamma$$
which is up to homotopy characterized by the property that it induces the identity on homology.
The symmetric Poincar{\'e} structure on 
$C_*(\overline{X};\rationals) \otimes_{\rationals \Gamma}
\mathcal{U}\Gamma$ induces one on
$ H_*(C_*(\overline{X};\rationals) \otimes_{\rationals \Gamma} \mathcal{U}\Gamma)$ 
and $i_*$ is an $\mathcal{U}\Gamma$-chain homotopy equivalence 
of symmetric  $\mathcal{U}\Gamma$-Poincar{\'e} complexes.
This implies for their classes in $L^0(\mathcal{U})$
\cite[Proposition 1.2.1]{Ranicki(1981)}.
$$[C_*(\overline{X};\rationals) \otimes_{\rationals \Gamma} \mathcal{U}\Gamma] ~ = ~
[H_*(C_*(\overline{X};\rationals) \otimes_{\rationals \Gamma} \mathcal{U}\Gamma)].$$
Elementary algebraic surgery in the sense of
\cite[Section 1.5]{Ranicki(1981)} shows that the class
$[H_*(C_*(\overline{X};\rationals) 
\otimes_{\rationals \Gamma} \mathcal{U}\Gamma)]$ in $L^0(\mathcal{U})$
is given by the sesquilinear non-degenerate pairing on the middle homology group
$H_{2n}(C_*(\overline{X};\rationals) \otimes_{\rationals \Gamma} \mathcal{U}\Gamma)$. 
Let 
\begin{equation*}
{\bf P}H_{2n}(C_*(\overline{X};\rationals) 
\otimes_{\rationals \Gamma} \NeumannN\Gamma)
\end{equation*}
be the projective part of
the finitely generated $\NeumannN \Gamma$-module 
$H_{2n}(C_*(\overline{X};\rationals) \otimes_{\rationals \Gamma} \NeumannN\Gamma)$ in the sense of
\cite[Definition 6.1]{Lueck(2002)}. It is a finitely generated projective
$\NeumannN (\Gamma)$-module \cite[Theorem 6.7]{Lueck(2002)}
and inherits a sesquilinear non-degenerate pairing from the Poincar{\'e} structure.
There is a canonical isomorphism
$$\left({\bf P} H_{2n}(C_*(\overline{X};\rationals) 
\otimes_{\rationals \Gamma} \NeumannN\Gamma)\right)
\otimes_{\NeumannN\Gamma} \mathcal{U}\Gamma)
\xrightarrow{\cong} H_{2n}(C_*(\overline{X};\rationals) 
\otimes_{\rationals \Gamma} \mathcal{U}\Gamma)$$
which is compatible with the pairings (see \cite[Theorem 6.7 and Lemma 8.33]{Lueck(2002)}). 
We have shown that the image of $\sigma(\overline{X})$ under the change of rings maps
$L^0(\rationals \Gamma) \to L^0(\NeumannN \Gamma)$ agrees with the class
represented by the Poincar{\'e} pairing on 
\begin{equation*}
{\bf P}H_{2n}(C_*(\overline{X};\rationals) \otimes_{\rationals \Gamma}
\NeumannN\Gamma).
\end{equation*}
We conclude from the definitions, Proposition \ref{PD_map}  and \cite[Theorem 6.24]{Lueck(2002)} 
that the map 
$$L^0(\NeumannN \Gamma) \xrightarrow{\cong} K_0(\NeumannN \Gamma) \to \reals$$
sends the class
represented by the Poincar{\'e} pairing on 
${\bf P}H_{2n}(C_*(\overline{X};\rationals) \otimes_{\rationals \Gamma} \NeumannN\Gamma)$ 
to $\sign^{(2)}_{\schain}(\overline{X})$. 
We conclude from the definition of $\sign_L^{(2)}(\overline{X})$
that $\sign^{(2)}_{\schain}(\overline{X}) = \sign_L^{(2)}(\overline{X})$ holds.
\end{proof}

If $X$ is a closed oriented smooth Riemannian
manifold then, as we have seen above, the signature operator
twisted with the canonical non-trivial flat $\NeumannN\Gamma$-bundle on $X$
has an index in $K_0(\NeumannN\Gamma)$. 

It is now a fundamental result, due to Mishchenko and Kasparov, that
this index is equal to the element given by the symmetric
signature (they are actually using the group $C^*$-algebra
$C^*\Gamma$, but the argument for the von Neumann algebra is the
same). For an extensive treatment of these facts (and a
generalization to more general $C^*$-algebra-module bundles), compare
\cite{Miller(1998)}. In particular, we get the following result
(see also \cite[pages 728-729]{Kreck-Leichtnam-Lueck(2002)}).

\begin{theorem}\label{theo:symmetric_signature_equals_K_index}
  Let $M$ be a closed oriented smooth Riemannian manifold of dimension $4n$.
  Let $\overline{M} \to M$ be a regular $\Gamma$-covering. Then
  \begin{equation*}
  \sign_L^{(2)}(\overline{M}) = \sign_K^{(2)}(\overline M).
  \end{equation*}
  \end{theorem}


\subsection{K\protect{\"u}nneth formula}
\label{sec:Kuenneth}

\begin{proposition}\label{sigmult}
  The $L^2$-signature is multiplicative: if $X$ and $Y$ are two
  Poincar{\'e} spaces with a regular $\Gamma_X$-covering $\overline{X} \to X$ and 
  a regular $\Gamma_Y$-covering $\overline{Y} \to Y$,
  then we get a regular $\Gamma_X \times \Gamma_Y$-covering 
  $\overline{X} \times \overline{Y} \to X \times Y$ and we have
  \begin{equation*}
    \sign^{(2)}(\overline{X}\times \overline{Y}) = \sign^{(2)}(\overline{X})\cdot
    \sign^{(2)}(\overline{Y}).
  \end{equation*}
\end{proposition}
\begin{proof}
  This follows, as in the classical compact case, in a straightforward
  way from the
  K{\"u}nneth formula for $L^2$-cohomology.
\end{proof}

\bibliographystyle{plain}
\bibliography{topol_L2_sig}

\def\cprime{$'$}
\begin{thebibliography}{10}

\bibitem{Aravinda-Farrell-Roushon(1997)}
C.S. Aravinda, F.T. Farrell, and S.K. Roushon.
\newblock Surgery groups of knot and link complements.
\newblock {\em Bull. of the London Math. Soc.}, 29:400--406, 1997.

\bibitem{Atiyah(1976)}
M.~F. Atiyah.
\newblock Elliptic operators, discrete groups and von {N}eumann algebras.
\newblock In {\em Colloque ``Analyse et Topologie'' en l'Honneur de Henri
  Cartan (Orsay, 1974)}, pages 43--72. Ast\'erisque, No. 32--33. Soc. Math.
  France, Paris, 1976.

\bibitem{Atiyah-Patodi-Singer(1975a)}
M.~F. Atiyah, V.~K. Patodi, and I.~M. Singer.
\newblock Spectral asymmetry and {R}iemannian geometry. {I}.
\newblock {\em Math. Proc. Cambridge Philos. Soc.}, 77:43--69, 1975.

\bibitem{Atiyah-Patodi-Singer(1975b)}
M.~F. Atiyah, V.~K. Patodi, and I.~M. Singer.
\newblock Spectral asymmetry and {R}iemannian geometry. {I}{I}.
\newblock {\em Math. Proc. Cambridge Philos. Soc.}, 78(3):405--432, 1975.

\bibitem{Atiyah-Patodi-Singer(1975c)}
M.~F. Atiyah, V.~K. Patodi, and I.~M. Singer.
\newblock Spectral asymmetry and {R}iemannian geometry. {I}{I}{I}.
\newblock {\em Math. Proc. Cambridge Philos. Soc.}, 79(1):71--99, 1976.

\bibitem{Dodziuk(1977)}
J.~Dodziuk.
\newblock De {R}ham-hodge theory for ${L}^2$-cohomology of infinite coverings.
\newblock {\em Topology}, 16:157--165, 1977.

\bibitem{Dold(1968)}
Albrecht Dold.
\newblock {\em On general cohomology. {C}hapters 1--9}.
\newblock Matematisk Institut, Aarhus Universitet, Aarhus, 1968.

\bibitem{Dupont(1978)}
J.~L. Dupont.
\newblock {\em Curvature and characteristic classes}.
\newblock Springer-Verlag, Berlin, 1978.
\newblock Lecture Notes in Mathematics, Vol. 640.

\bibitem{Farrell-Jones(1988b)}
F.T. Farrell and L.E. Jones.
\newblock The surgery {L}-groups of poly-(finite or cyclic) groups.
\newblock {\em Inventiones Mathematicae}, 91:559--586, 1988.

\bibitem{Farrell-Jones(1993c)}
F.T. Farrell and L.E. Jones.
\newblock Topological rigidity for compact non-positively curved manifolds.
\newblock In Robert Greene, editor, {\em Differential geometry. Part 3:
  Riemannian geometry. Proceedings of a summer research institute, held at the
  UCLA, July 8-28, 1990}, volume 54 part 3 of {\em Proc. of Symp. in Pure
  Math.}, pages 229--274. AMS, 1993.
\newblock Zbl. 796.53043.

\bibitem{Gaffney(1954)}
M.~P. Gaffney.
\newblock A special {S}tokes's theorem for complete {R}iemannian manifolds.
\newblock {\em Ann. of Math. (2)}, 60:140--145, 1954.

\bibitem{Goldstein_et_al(1988)}
V.~M. Gol{\cprime}dshte{\u\i}n, V.~I. Kuz{\cprime}minov, and I.~A. Shvedov.
\newblock The de {R}ham isomorphism of the ${L}\sb p$-cohomology of noncompact
  {R}iemannian manifolds.
\newblock {\em Sibirsk. Mat. Zh.}, 29(2):34--44, 216, 1988.

\bibitem{Higson-Kasparov(2001)}
N.~Higson and G.~Kasparov.
\newblock ${E}$-theory and ${K}{K}$-theory for groups which act properly and
  isometrically on {H}ilbert space.
\newblock {\em Invent. Math.}, 144(1):23--74, 2001.

\bibitem{Hilsum(1985)}
M.~Hilsum.
\newblock Signature operator on {L}ipschitz manifolds and unbounded {K}asparov
  bimodules.
\newblock In {\em Operator algebras and their connections with topology and
  ergodic theory (Bu\c steni, 1983)}, pages 254--288. Springer, Berlin, 1985.

\bibitem{Hilsum(1989)}
M.~Hilsum.
\newblock Fonctorialit\'e en ${K}$-th\'eorie bivariante pour les vari\'et\'es
  lipschitziennes.
\newblock {\em K-theory}, 3:401--440, 1989.

\bibitem{Kreck-Leichtnam-Lueck(2002)}
M.~Kreck, E.~Leichtnam, and W.~L{\"u}ck.
\newblock On the cut and paste property of higher signatures of a closed
  oriented manifold.
\newblock {\em Topology}, 41:725--744, 2002.

\bibitem{Lott-Lueck(1995)}
J.~Lott and W.~L{\"u}ck.
\newblock $l^2$-topological invariants of $3$-manifolds.
\newblock {\em Inventiones Mathematicae}, 120:15--60, 1995.

\bibitem{Lueck(2002)}
W.~L\"uck.
\newblock {\em $L^2$-Invariants: Theory and Applications to Geometry and
  {K}-Theory}.
\newblock Springer-Verlag, Berlin, 2002.
\newblock Ergebnisse der Mathematik und ihrer Grenzgebiete Vol 44.

\bibitem{Lueck-Schick(1998)}
W.~L{\"u}ck and T.~Schick.
\newblock ${L}\sp 2$-torsion of hyperbolic manifolds of finite volume.
\newblock {\em Geom. Funct. Anal.}, 9(3):518--567, 1999.

\bibitem{Lueck-Schick(2001)}
W.~L\"uck and T.~Schick.
\newblock Approximating ${L}^2$-signatures by their compact analogues.
\newblock preprint, SFB Geometrische Strukturen, M\"unster,
  http://front.math.ucdavis.edu/math.GT/0110328, 2001.

\bibitem{Miller(1998)}
J.~G. Miller.
\newblock Signature operators and surgery groups over ${C}\sp *$-algebras.
\newblock {\em $K$-Theory}, 13(4):363--402, 1998.

\bibitem{Mishchenko(1970)}
A.~S. Mi{\v{s}}{\v{c}}enko.
\newblock Homotopy invariants of multiply connected manifolds. {I}. {R}ational
  invariants.
\newblock {\em Izv. Akad. Nauk SSSR Ser. Mat.}, 34:501--514, 1970.

\bibitem{Mislin(2002)}
G.~Mislin.
\newblock Equivariant {K}-homology of the classifying space for proper actions.
\newblock Unpublished preprint, manuscript for a school in Barcelona held in
  September 2001.

\bibitem{Moscovici-Wu(1993)}
H.~Moscovici and F.~Wu.
\newblock Localization of topological {P}ontrjagin classes via finite
  propagation speed.
\newblock {\em C.R. Acad. Sci. Paris, s\'eries 1}, 317(7):661--665, 1993.

\bibitem{Ramachandran(1993)}
M.~Ramachandran.
\newblock von {N}eumann index theorems for manifolds with boundary.
\newblock {\em J. Differential Geom.}, 38(2):315--349, 1993.

\bibitem{Ranicki(1980aa)}
A.~Ranicki.
\newblock The algebraic theory of surgery. {I}. {F}oundations.
\newblock {\em Proc. London Math. Soc. (3)}, 40(1):87--192, 1980.

\bibitem{Ranicki(1980bb)}
A.~Ranicki.
\newblock The algebraic theory of surgery. {I}{I}. {A}pplications to topology.
\newblock {\em Proc. London Math. Soc. (3)}, 40(2):193--283, 1980.

\bibitem{Ranicki(1981)}
A.~Ranicki.
\newblock {\em Exact sequences in the algebraic theory of surgery}.
\newblock Princeton University Press, 1981.

\bibitem{Ranicki(1992)}
A.~Ranicki.
\newblock {\em Algebraic ${L}$-theory and topological manifolds}, volume 102 of
  {\em Cambridge Tracts in Mathematics}.
\newblock Cambridge University Press, 1992.

\bibitem{Ranicki(1995b)}
A.~Ranicki.
\newblock On the {N}ovikov conjecture.
\newblock In {\em Novikov conjectures, index theorems and rigidity, Vol.\ 1
  ({O}berwolfach, 1993)}, pages 272--337. Cambridge Univ. Press, Cambridge,
  1995.

\bibitem{Ranicki-Sullivan(1976)}
A.~Ranicki and D.~Sullivan.
\newblock A semi-local combinatorial formula for the signature of a
  $4k$-manifold.
\newblock {\em J. Differential Geometry}, 11(1):23--29, 1976.

\bibitem{Rinow(1975)}
W.~Rinow.
\newblock {\em Lehrbuch der {T}opologie}.
\newblock VEB Deutscher Verlag der Wissenschaften, Berlin, 1975.
\newblock Hochschulb\"ucher f\"ur Mathematik, Band 79.

\bibitem{Rosenberg(1983)}
J.~Rosenberg.
\newblock ${C}\sp{\ast} $-algebras, positive scalar curvature, and the
  {N}ovikov conjecture.
\newblock {\em Inst. Hautes \'Etudes Sci. Publ. Math.}, 58:197--212 (1984),
  1983.

\bibitem{Rosenberg(1995)}
J.~Rosenberg.
\newblock Analytic {N}ovikov for topologists.
\newblock In {\em Proceedings of the conference ``Novikov conjectures, index
  theorems and rigidity'' volume I, Oberwolfach 1993}, volume 226 of {\em LMS
  Lecture Notes Series}, pages 338--372. Cambridge University Press, 1995.

\bibitem{Schafer(1970)}
J.~A. Schafer.
\newblock Topological {P}ontrjagin classes.
\newblock {\em Comment. Math. Helv.}, 45:315--332, 1970.

\bibitem{Schick(2002)}
T.~Schick.
\newblock A {KK}-proof of {A}tiyah's ${L}^2$-index theorem.
\newblock 2002, in preparation.

\bibitem{Schick(1996)url}
T.~Schick.
\newblock {\em Analysis on $\partial$-manifolds of bounded geometry, {H}odge-de
  {R}ham isomorphism and $L^2$-index theorem}.
\newblock Shaker, Aachen, 1996.
\newblock (Dissertation, Mainz),
  http://www.uni-math.gwdg.de/schick/publ/dissschick.html.

\bibitem{Schick(2001)}
T.~Schick.
\newblock ${L}\sp 2$-index theorem for elliptic differential boundary problems.
\newblock {\em Pacific J. Math.}, 197(2):423--439, 2001.

\bibitem{Sullivan(1979)}
D.~Sullivan.
\newblock Hyperbolic geometry and homeomorphisms.
\newblock In {\em Geometric topology (Proc. Georgia Topology Conf., Athens,
  Ga., 1977)}, pages 543--555. Academic Press, New York, 1979.

\bibitem{Teleman(1983)}
N.~Teleman.
\newblock The index of signature operators on {L}ipschitz manifolds.
\newblock {\em Publ. Math. IHES}, 58:39--78, 1983.

\bibitem{Teleman(1984)}
N.~Teleman.
\newblock The index theorem for topological manifolds.
\newblock {\em Acta Math.}, 153(1-2):117--152, 1984.

\bibitem{Vaillant(1997)}
B.~Vaillant.
\newblock Indextheorie f\"ur {{\"U}}berlagerungen.
\newblock Diplomarbeit, Universit\"at Bonn,
  http://styx.math.uni-bonn.de/boris/diplom.html, 1997.

\bibitem{Valette(2002)}
A.~Valette.
\newblock {\em Introduction to the {B}aum-{C}onnes conjecture}.
\newblock Birkh\"auser Verlag, Basel, 2002.
\newblock From notes taken by Indira Chatterji, With an appendix by Guido
  Mislin.

\bibitem{Wall(1966)}
C.~T.~C. Wall.
\newblock Surgery of non-simply-connected manifolds.
\newblock {\em Ann. of Math. (2)}, 84:217--276, 1966.

\bibitem{Wall(1967)}
C.~T.~C. Wall.
\newblock Poincar\'e complexes. {I}.
\newblock {\em Ann. of Math. (2)}, 86:213--245, 1967.

\bibitem{Wall(1999)}
C.~T.~C. Wall.
\newblock {\em Surgery on compact manifolds}.
\newblock AMS, second edition, 1999.
\newblock Edited and with a foreword by A. A. Ranicki.

\bibitem{Weinberger(1988b)}
S.~Weinberger.
\newblock Homotopy invariance of $\eta$-invariants.
\newblock {\em Proc. Nat. Acad. Sci.}, 85:5362--5363, 1988.

\end{thebibliography}

-----------------------------7864765581673488648746260076--

\end{document}